\documentclass[twoside, 12pt]{article}
\usepackage[T1]{fontenc}
\usepackage[latin1]{inputenc}
\usepackage{amscd}
\usepackage{amsmath}
\usepackage{amsfonts}
\usepackage{amssymb}
\usepackage{amsthm}
\usepackage{comment}
\usepackage{graphicx,epsfig}
\usepackage{changebar}

\renewcommand{\tilde}{\widetilde}
\renewcommand{\hat}{\widehat}
\newcommand{\moins}{\backslash}

\newcommand{\norm}[1]{\left\| #1 \right\|}

    \DeclareMathOperator{\Card}{Card}

    \DeclareMathOperator{\Leb}{Leb}
    
    \DeclareMathOperator{\dist}{dist}

\newcommand{\Sing}{Sing}
\newcommand{\de}{{\rm d}}
\newcommand{\dd}{\, {\rm d}}

\renewcommand{\geq}{\geqslant}
\renewcommand{\leq}{\leqslant}

\newcommand{\N}{\mathbb{N}}
\newcommand{\Z}{\mathbb{Z}}
\newcommand{\R}{\mathbb{R}}

\newcommand{\C}{\mathbb{C}}

\renewcommand{\phi}{\varphi}
\renewcommand{\epsilon}{\varepsilon}
\newcommand{\tq}{\ |\ }

\newcommand{\boC}{\mathcal{C}}

\newcommand{\boD}{\mathcal{D}}
\newcommand{\boB}{\mathcal{B}}
\newcommand{\boN}{\mathcal{N}}

\newcommand{\boH}{\mathcal{H}}

\newtheorem{thm}{Theorem}[section]
\newtheorem*{thm*}{Theorem}
\newtheorem{prop}[thm]{Proposition}
\newtheorem{lem}[thm]{Lemma}
\newtheorem{cor}[thm]{Corollary}
\newtheorem{defn}[thm]{Definition}

\theoremstyle{definition}

\newtheorem{rmq}[thm]{Remark}
\newtheorem*{rmq*}{Remark}

\newcommand{\Ca}{C}
\newcommand{\Cb}{C}

\newcommand{\Cd}{C}
\newcommand{\Ce}{C_3}
\newcommand{\Cf}{C_4}
\newcommand{\Cg}{C_5}
\newcommand{\Ch}{C}

\newcommand{\Cz}{C}
\newcommand{\Cy}{C}
\newcommand{\Cx}{C}
\newcommand{\Cw}{C}
\newcommand{\Cv}{C}

\newcommand{\Cmm}{C_1}
\newcommand{\Cn}{C_2}
\newcommand{\Co}{C}

\newcommand{\rhp}{s}
\newcommand{\gam}{h}
\newcommand{\Gam}{H}
\newcommand{\haut}{\omega}

\title{Limit theorems in the stadium billiard}

\author{P\'eter B\'alint\thanks{
\emph{Address}: Institute of Mathematics, Budapest University of
Technology and Economics, H-1111 Egry J\'ozsef u. 1, Budapest,
Hungary; \emph{Email:} bp@renyi.hu; \emph{URL:}
http://www.renyi.hu/{\lower.7ex\hbox{\~{}}}bp/} \  and S\'ebastien
Gou\"ezel\thanks{ \emph{Address:} D\'epartement de Math\'ematiques et
Applications, École Normale Sup\'erieure, 45 rue d'Ulm, Paris,
France; \emph{Email:} Sebastien.Gouezel@ens.fr; \emph{URL:}
http://www.dma.ens.fr/{\lower.7ex\hbox{\~{}}}gouezel/} }

\begin{document}

\maketitle

\begin{abstract}
We prove that the Birkhoff sums for ``almost every'' relevant
observable in the stadium billiard obey a non-standard limit law.
More precisely, the usual central limit theorem holds for an
observable if and only if its integral along a one-codimensional
invariant set vanishes, otherwise a $\sqrt{n\log n}$ normalization
is needed. As one of the two key steps in the argument, we obtain
a limit theorem that holds in Young towers with exponential return
time statistics in general, an abstract result that seems to be
applicable to many other situations.
\end{abstract}

\section*{Introduction}

The subject of this article, the stadium billiard, belongs to the
class of dynamical systems that are sometimes referred to as
intermittent ones. This name is related to the \emph{weakly
chaotic} nature of the time evolution that accounts for a
modified, relaxed appearance of the behavior characteristic to
systems with uniform hyperbolicity. In particular, the
mathematically rigorous investigation of the stadium started with
\cite{bunimovich:stadium} where Bunimovich showed (with respect to
the natural invariant measure) that the Lyapunov exponents are
almost everywhere non-zero, and that the system is ergodic. Thus
in that respect the stadium billiard resembles dispersing
billiards, however, when finer statistical properties are
discussed, deviations start to show up. Recent works by Markarian
(\cite{markarian:slow}) and Chernov--Zhang (\cite{chernov:slow})
have obtained an upper bound on the rate of mixing: given two
sufficiently smooth (H\"older or Lipschitz continuous) observables,
their correlations decay as $O((\log n)^2/n)$. Although this upper
bound is most likely not sharp, it is definitely not far from the
optimal either (see Corollary~\ref{cor:lower_bound_corr}). In this
paper we investigate the issue of probabilistic limit laws and
provide further evidence of the intermittent nature of the
dynamics. Namely we show that the limit behavior of a sufficiently
smooth observable with zero mean, to be denoted by $f_0$, is
characterized by a quantity $I$ (cf.\ \eqref{donne_I}), its
average along the one dimensional set of trajectories bouncing
forever along the straight segments. In the typical case $I\ne 0$,
the Birkhoff sums of $f_0$ satisfy a non-standard limit theorem --
convergence in distribution to the Gaussian law can be obtained
with a $\sqrt{cn\log n}$ normalization, where the constant $c$ is
a multiple of $I^2$, see Theorem~\ref{main_theorem}. On the other
hand the central limit theorem in its usual form applies if $I=0$,
see Theorem~\ref{main_thm_2}. These results have some almost
immediate corollaries: we obtain the analogous limit theorems for
the billiard flow (Corollary~\ref{cor:flow}) and, though in a very
weak form, some lower bounds on the rate of correlation decay
(Corollary~\ref{cor:lower_bound_corr}).

The issue of probabilistic limit laws in dynamical systems has a
long history. In the chaotic setting the possibly most frequently
applied method is Gordin's martingale argument (see \cite{gordin},
or \cite{lsyoung:annals} and references therein) that roughly
states that under quite general conditions, whenever the
correlations decay at a summable rate, the usual central limit
theorem holds. This technique, however, cannot treat non-standard
limit behavior or non-summable decay rates. Recently Aaronson and
Denker have proposed an approach to the issue of non-standard
limit theorems, see e.g.\ \cite{aaronson_denker:central}. The
dynamical systems they study, the so called Gibbs-Markov maps,
possess some important features characteristic to uniformly
expanding Markov maps of the interval, in particular, they are
strongly chaotic. However, the functions $f$ for which limit
theorems are proved are unbounded, and do not even belong to
$L^2$. This setting allows for the use of Perron-Frobenius
techniques: there is a one parameter family of transfer operators
the spectra of which give precise information on the limit
behavior of the observable. In particular, the Birkhoff sums
satisfy exactly the same limit theorem that an i.i.d.\ sequence of
random variables with the distribution of $f$ would have. For
details see \cite{aaronson_denker:central} and
Section~\ref{par_attraction} of the present paper.

The above ideas can be implemented to treat limit laws for bounded
functions in weakly chaotic systems $T_0:X_0\to X_0$ in case the
following scenario applies. Let us assume that the source of
non-uniformity in hyperbolicity is a well-distinguishable
geometric effect. Then one may consider a subset $X\subset X_0$
such that the first return map onto $X$ is uniformly hyperbolic,
however, our observable induces an unbounded function on $X$. Thus
we arrive at a setting close to that of
\cite{aaronson_denker:central}. This line of approach has been
successfully applied to systems for which the induced map is
Gibbs-Markov (see eg. \cite{gouezel:stable}), which, however, is
not exactly the case of the stadium billiard.

What replaces Gibbs-Markov property in billiards is the presence
of a Young tower, an object that has turned out to be very
effective when estimating the rate of the decay of correlations.
There are two versions of Young towers: those with exponential
return time statistics ensure rapid mixing -- exponential decay of
correlations -- via Perron-Frobenius techniques
(\cite{lsyoung:annals}), while those with polynomial return time
statistics give polynomial upper bounds on the rate of correlation
decay -- slow mixing rates -- via coupling techniques
(\cite{lsyoung:recurrence}). As to the case of the stadium
billiard, the Young towers constructed in \cite{markarian:slow}
and \cite{chernov:slow} have polynomial return time statistics
with respect to the original map, and exponential return time
statistics with respect to the induced map. The aim of the present
paper is, in addition to present our results on the stadium
billiard, to demonstrate that Young towers, originally designed to
estimate mixing rates, are almost equally powerful when the issue
of various limit laws is investigated. Note that this fact has
already been observed and emphasized by Sz\'asz and Varj\'u in the
papers \cite{szasz_varju:finite} and
\cite{szasz_varju:expository}.

The proof of Theorem~\ref{main_theorem} consists of two clearly
distinguishable ingredients. On the one hand, via Perron-Frobenius
techniques, we prove Theorem~\ref{thm_ad_young}, a general result
in Young towers with exponential return time statistics. This
concerns the limit behavior of the Birkhoff sums of observables
belonging to the non-standard domain of attraction of the Gaussian
law. It is important to note that, as the Gibbs-Markov property is
replaced by a Young tower, a new effect shows up that typically
rescales the normalizing sequence with a constant multiplicator.
We would also like to emphasize that this first ingredient of the
proof is completely general and could be applied to many other
situations. On the other hand, the second ingredient is directly
related to the stadium billiard. We rely on \cite{markarian:slow}
and \cite{chernov:slow} when considering a suitable induced map
that allows for a Young tower with exponential return time
statistics. However, in order to ``pull back'' the limit theorem
from the Young tower to the phase space of the billiard, and in
order to give a transparent interpretation in terms of quantities
easy to calculate, we need to perform a finer and more detailed
geometric analysis of the stadium than the one presented in the
above two papers.

We strongly do believe that our line of approach could be applied
to obtain non-standard limit theorems in many other hyperbolic
dynamical systems, in particular, in certain billiards with slow
mixing rates. One of the most interesting candidates, the infinite
horizon Lorentz process, for which the significance of the limit
behavior is further emphasized as it may give an effective tool to
discuss recurrence properties, is investigated by Sz\'asz and
Varj\'u (\cite{szasz_varju:infinite}). Among others, it is also
worth mentioning skewed stadia (see \cite{chernov:slow}) and
dispersing billiards with cusps (\cite{machta:cusp}). We plan to
turn back to these systems in separate papers.

The article has five sections. In the first one we state our main
results and fix some basic notation. The second section is devoted
to general results on the stadium billiard. We essentially recall
the existence of Young towers for an induced map, proved by
Markarian in \cite{markarian:slow}. In the third part, we study
abstract Young towers and establish a spectral perturbation
estimate. In particular, to get a limit theorem, it is sufficient
to study an integral with sufficient precision. In Section~4, we
come back to the stadium billiard map, and describe geometrically
this integral. With a careful study of the singularities of the
stadium map, this gives an accurate description of this integral.
Finally, in Section~5, we use together the abstract results of
Section~3 and the explicit estimate of Section~4, to prove
Theorem~\ref{main_theorem}.

\section{Results}

Let $\ell>0$. We consider a region in the plane delimited by two
semicircles of radius $1$, joined by two horizontal segments of
length $\ell$, tangent to the semicircles. 
To a point on the boundary of this set and a vector
pointing inwards, we associate an image by the usual billiard
reflection law. This defines the stadium billiard map $T_0 : X_0
\to X_0$. This map admits a unique absolutely continuous invariant
probability measure $\mu_0$.

A point in the phase space $X_0$ is given by $(r,\theta)$, where
$r\in \R/(2\pi+2\ell)\Z$ is the position on the boundary, and
$\theta \in (-\pi/2,\pi/2)$ is the angle with respect to the
normal to this boundary at $r$. The invariant measure $\mu_0$ is
given by
  \begin{equation*}
  \dd\mu_0=\frac{\cos \theta \dd r \dd \theta}{2(2\pi+2\ell)}.
  \end{equation*}

We will assume that $r=0$ corresponds to the lower endpoint of the
right semi-circle, and that the boundary is oriented
counterclockwise. Hence, the semicircles correspond to $0\leq r
\leq \pi$ and $\pi+\ell \leq r \leq 2\pi+\ell$.

Let $f_0:X_0 \to \R$ be a H\"older function. We will be interested
in the asymptotic behavior of the Birkhoff sums of $f_0$. The map
$T_0$ is slowly mixing, by \cite{markarian:slow} and
\cite{chernov:slow}: its correlations decay (at least) like
$O((\log n)^2/n)$. This estimate is not summable, whence the usual
Gordin martingale argument to get a central limit theorem does not
apply. We will indeed prove that the usual central limit theorem
does not hold.

Let
  \begin{equation}
  \label{donne_I}
  I=\frac{1}{2\ell}\left[
  \int_{r\in[\pi,\pi+\ell]} f_0(r,0)\dd r  +
  \int_{r\in [2\pi+\ell,2\pi+2\ell]} f_0(r,0)\dd r\right].
  \end{equation}
This is the average of $f_0$ along the trajectories bouncing
perpendicularly to the segments of the stadium.

In this article, we prove the following theorem:
\begin{thm}
\label{main_theorem}
Let $f_0:X_0 \to \R$ be H\"older continuous, satisfying $\int f_0
\dd\mu_0=0$ and $I\not=0$. Then
  \begin{equation*}
  \frac{\sum_{k=0}^{n-1} f_0\circ T_0^k}
   {\sqrt{c n \log n }} \to \boN(0,1),
  \end{equation*}
where
  \begin{equation*}
  c= \frac{ 4+3\log 3}{4-3\log 3} \cdot \frac{ \ell^2 I^2}{ 4
  (\pi+\ell)}.
  \end{equation*}
\end{thm}

\begin{cor}
\label{cor:lower_bound_corr}
Under the assumptions of Theorem~\ref{main_theorem}, the quantity
$n\int f_0 \cdot f_0 \circ T_0^n$ does not tend to zero.
\end{cor}
\begin{proof}
We have
  \begin{equation*}
  \int \left[ \sum_{k=0}^{n-1} f_0 \circ T_0^k \right]^2
  = n \int f_0^2+2\sum_{i=1}^{n-1} (n-i) \int f_0 \cdot f_0 \circ T^i.
  \end{equation*}
If $\int f_0 \cdot f_0 \circ T^i=o(1/i)$, we obtain $\int \left[
\sum_{k=0}^{n-1} f_0 \circ T_0^k \right]^2 = o(n \log n)$. In
particular, the variance of the random variable
$\frac{\sum_{k=0}^{n-1} f_0\circ T_0^k}{\sqrt{n \log n}}$ tends to
zero. This implies that this random variable tends to zero in
probability, which is in contradiction with
Theorem~\ref{main_theorem}.
\end{proof}

Hence, we obtain a lower bound $O(1/n)$ on the speed of decay of
correlations of H\"older functions.  It indicates that the upper
bound of Markarian and Chernov-Zhang is close to optimal (it may
probably be replaced by $O(1/n)$, since the $(\log n)^2$ seems to
be due to the technique of proof).

We also obtain the following (easier) result:
\begin{thm}
\label{main_thm_2}
Let $f_0:X_0 \to \R$ be H\"older continuous, satisfying $\int f_0
\dd\mu_0=0$ and $I=0$. Then there exists $\sigma^2\geq 0$ such
that
  \begin{equation*}
  \frac{\sum_{k=0}^{n-1} f_0\circ T_0^k}
   {\sqrt{n}} \to \boN(0,\sigma^2).
  \end{equation*}
\end{thm}
Hence, when $I=0$, the Birkhoff sums of $f_0$ satisfy a usual
central limit theorem.

Before going into the details of the proof we consider one
particularly interesting observable: the free path. Given
$x=(r,\theta)$, we denote $T_0x=(r_1,\theta_1)$ and define
$\tau(x)$ as the planar distance of $r$ and $r_1$. In other words,
$\tau(x)$ is the length of the trajectory segment the point
particle follows until the next collision. To investigate the
limit behavior of the free path $\tau: X_0\to \R$, we have to
subtract its mean
${\bar \tau}=\int \tau \dd\mu_0$, thus we
define $\tau_0(x)=\tau(x) -{\bar \tau}$. There is a remarkably
simple formula for ${\bar \tau}$ that can be obtained by comparing
the invariant measures for the billiard map and the billiard flow
(see \cite{chernov:entropy}):
\begin{equation}
\label{freepath}
{\bar \tau}= \frac{\pi(\pi+2\ell)}{2\ell+2\pi}.
\end{equation}
On the other hand, we may easily calculate (\ref{donne_I}) as we
have $\tau(r,0)=2$ whenever $r\in[\pi,\pi+\ell]$ or $r\in
[2\pi+\ell,2\pi+2\ell]$, thus $I_{\tau}=2$ and
$I_{\tau_0}=2-\frac{\pi(\pi+2\ell)}{2\ell+2\pi}$. This means there
is a ``best'' stadium with
$\ell=\ell^*=\frac{4\pi-\pi^2}{2\pi-4}\approx 1.18$ for which
$I_{\tau_0}=0$ and consequently, by Theorem~\ref{main_thm_2} the
(centralized) free path satisfies the usual central limit theorem.
However, whenever $\ell\ne \ell^*$, we have $I_{\tau_0}\ne 0$ and,
by Theorem~\ref{main_theorem} a stronger normalization is needed.

Our interest in $\tau$ is also related to the fact that the
billiard flow may be considered as a suspension above the billiard
map with the roof function $\tau(x)$. By \cite{melbourne_torok}
suspension flows do inherit some statistical properties 
from the base transformation, in
particular limit theorems, under
quite general conditions. Let us denote the billiard flow by
\begin{equation*}
X_{\tau}=\{ (x,u) \tq x\in X_0, 0\le u\le \tau(x) \}/\sim, \quad\!
(x,\tau(x))\sim (T_0x,0) \quad S_t(x,u)=(x,u+t), \quad\!
\mu_{\tau}=\mu_0\times\frac{\Leb}{{\bar \tau}}
\end{equation*}
where the action of the flow is understood modulo identifications.
Consider a H\"older observable $\Phi:X_{\tau}\to \R$ satisfying
$\int \Phi \dd\mu_{\tau}=0$, and define
\begin{equation*}
\Phi_T(x) =\int_0^T \Phi(S^tx)\dd t; \qquad J_{\Phi}=
\frac{1}{4\ell}\left[
  \int_{r\in[\pi,\pi+\ell]\cup[2\pi+\ell,2\pi+2\ell]}
  \int_{ t\in [0,2]} \Phi(r,0,t)\dd t\dd r  \right].
\end{equation*}

\begin{cor}
\label{cor:flow}
\begin{enumerate}
\item If $J_{\Phi}\ne 0$, then
\begin{equation*}
\frac{\Phi_T}{\sqrt{\frac{c}{\bar \tau}T\log T}} \to \boN(0,1).
\end{equation*}
Here $c$ is the constant from Theorem~\ref{main_theorem}, with $I$
replaced by $J_{\Phi}$.
\item If $J_{\Phi}=0$, then
\begin{equation*}
  \frac{\Phi_T}
   {\sqrt{T}} \to \boN(0,\sigma_{\Phi}^2)
\end{equation*}
for some $\sigma_{\Phi}^2\geq 0$.
\end{enumerate}
\end{cor}

\begin{proof}
Define $f_0:X_0\to \R$ as $f_0(x)=\int_0^{\tau(x)} \Phi(S_t(x,0))
\dd t$. Then $f_0$ is H\"older, $\int f_0 \dd \mu_0=0$ and
$I_{f_0}=J_{\Phi}$. Thus, depending on the value of $J_{\Phi}$,
one of our two main theorems applies. To show that $\Phi$ inherits
the limit behavior from $f_0$, we apply the flow version of
\cite[Theorem A1]{gouezel:skewproduct} recalled as
Theorem~\ref{thm_probabiliste_general} in this paper (see also
Remark~\ref{rmq_limit_flows}). We only need to check that the
three conditions of this theorem are satisfied. In case
$J_{\Phi}\ne 0$ (and even if $J_{\Phi}=0$ and $\ell=\ell^*$)
conditions 1 and 3 are satisfied with $b=1$. Then condition 2 is
merely the Birkhoff ergodic theorem, thus the first statement is
established. If $J_{\Phi}=0$, the appropriate normalization for
$\tau$ may be $\sqrt{n\log n}$ as opposed to $\sqrt{n}$ needed for
$f_0$. Thus conditions 1 and 3 of
Theorem~\ref{thm_probabiliste_general} are satisfied for any
$0<b<1$, but not for $b=1$. This means condition 2 is to be
established for some $b<1$, but this is merely our
Remark~\ref{rmq:strong_ergthm_f0}. This completes the proof of the
second statement.
\end{proof}

We will say that a H\"older continuous function $f_0:X_0 \to \R$
with vanishing integral satisfies $(P1)$ if $I\not=0$ and $f_0$
vanishes on the set of points $x$ such that $x, T_0(x)$ and
$T_0^{-1}(x)$ belong to the same semicircle, and that $f_0$ satisfies
$(P2)$ if $I=0$. We will in fact prove Theorem~\ref{main_theorem}
for functions satisfying $(P1)$, and Theorem~\ref{main_thm_2} for
functions satisfying $(P2)$. This will imply
Theorem~\ref{main_theorem} in full generality. Namely, if $f_0$ is
H\"older continuous and satisfies $I\not=0$, then we may write it as
$f_0=f_1+f_2$ where $f_1$ satisfies $(P1)$ and $f_2$ satisfies
$(P2)$. By Theorem~\ref{main_thm_2}, $\frac{S_n f_2}{\sqrt{n \log
n}} \to 0$. Hence, it is equivalent to have
Theorem~\ref{main_theorem} for $f_0$ or $f_1$. We will 
comment on the technical reason for introducing the classes (P1) and (P2) 
in Remark~\ref{rmq:sliding_nocoboundary} below.

In this paper, $C$ will denote a generic constant, that can change
from one occurrence to the next. Some constants, which will be
used at different places in the paper, will be denoted by
$C_1,C_2,\dots$ and will have a fixed value.

\section{Background material on the  stadium billiard}

\label{section_background}

\subsection{Geometric description of the initial map and of an induced
map}
\label{sec:inigeom}

The map $T_0$ has almost everywhere two nonzero Lyapunov exponent.
However, the expansion in the unstable cone (and the contraction
in the stable cone) are not uniform: points bouncing many times
along the segments, or sliding along the circles, have an
expansion arbitrarily close to $1$.

To get uniform expansion, we follow \cite{markarian:slow} and
\cite{chernov:slow}. Let $X$ be the set of points $x$ in $X_0$
such that $x$ belongs to a semicircle and $T^{-1}(x)$ does not
belong to this semicircle. The set $X_0$ is a union of two
parallelograms in $X$, and it satisfies
  \begin{equation}
  \label{calcule_mu0_X}
  \mu_0(X)=\frac{2\pi}{2(2\pi+2\ell)}=\frac{\pi}{2(\pi+\ell)}.
  \end{equation}
Define a new probability measure on $X$ by
  \begin{equation*}
  \dd\mu= \frac{ \cos \theta \dd r \dd \theta}{2\pi}.
  \end{equation*}
For $x\in X$, let $\phi_+(x)=\inf\{ n\geq 1, T_0^n(x)\in X\}$.
This is the return time of $x$. Let $T:X \to X$ be the first
return map, induced by $T_0$ on $X$, i.e.,
$T(x)=T_0^{\phi_+(x)}(x)$. This map preserves the probability
measure $\mu$ on $X$. Moreover, it is uniformly hyperbolic in the
following sense:
\begin{prop}
\label{prop_hyperbolic}
There exists a continuous family of closed cones $\boC^u(x)$ for
$x\in X$, such that $DT(x)(\boC^u(x))\subset \boC^u(Tx)$.
Moreover, there exist constants $\Lambda>1$ and $C>0$ such that,
for all $x\in X$, for all $v\in \boC^u(x)$, for all $n\in \N$ such
that $T^n$ is defined and differentiable at $x$,
  \begin{equation*}
  \norm{DT^n(x)v} \geq C \Lambda^n \norm{v}.
  \end{equation*}
Moreover, these cones are uniformly bounded away from the
horizontal and vertical directions (i.e., $\{\de \theta=0\}$ and
$\{\de r=0\}$).

In the same way, there exist stable cones $\boC^s(x)$, which
satisfy the same properties for $T^{-1}$, except that they are not
bounded away from the horizontal direction.
\end{prop}
This proposition can be found in \cite{markarian:slow} and
\cite{chernov:slow}. The uniform expansion is given for a
different metric, the $p$-metric, satisfying $\norm{v}_p \leq
\norm{v}$. However, it is easy to check that there exists $C>0$
such that, for all $x\in X$, for all $v\in \boC^u(x)$,
$\norm{DT(x)v}_p \geq C \norm{v}$. Hence, the uniform expansion in
the $p$-metric implies the same statement for the euclidean
metric, up to a constant $C$.

There are two different types of points for which $\phi_+(x)$ can
be large: they correspond to points bouncing many times along the
segments, or sliding many times along the circles. We will need to
describe rather precisely the hyperbolic behavior of $T$ in
bouncing regions:
\begin{prop}
\label{prop_contracts}
If $x$ is a bouncing point satisfying $\phi_+(x)=n$, then $T$
contracts the p-metric of vectors in the stable cone at least by a
factor $\frac{C}{n}$, while $T^{-1}$ contracts the p-metric of
vectors in the unstable cone at least by a factor $\frac{C}{n}$.
Moreover, $Tx$ and $T^{-1}x$ are bouncing points with $\phi_+(Tx)
\geq n/4$, $\phi_+(T^{-1}x) \geq n/4$ if $n$ is large enough. This
implies, in turn, that the above contraction estimates are valid
in the euclidean metric as well.
\end{prop}

\begin{rmq}
\label{rmq:sliding_nocoboundary}
Note that if $x$ is a sliding point satisfying $\phi_+(x)=n$, then we 
can only guarrantee that $Tx$ and $T^{-1}x$ are sliding points with 
$\phi_+(Tx)\geq C \sqrt{n}$ and $\phi_+(T^{-1}x)\geq C \sqrt{n}$. 
This has an unfortunate consequence: we can only apply the coboundary 
arguments of section~\ref{sec:coboundary} to functions vanishing along sliding 
trajectories. Essentially this is the technical reason for introducing the 
classes (P1) and (P2). The proof of Theorem~\ref{main_theorem} relies 
heavily on Perron--Frobenius techniques, and thus requires an expanding 
setting, which implies that collapsing along stable manifolds -- 
coboundary arguments -- are essential. Now for the class (P2) it is enough 
to prove the usual central limit theorem (Theorem~\ref{main_thm_2}) which 
can be carried out in a roundabout way in the hyperbolic setting, see 
section~\ref{sec:gordin}.
\end{rmq}

\subsection{ Young tower of $T$}

A set $R \subset X$ is a \emph{rectangle} if there exist $x\in R$
with a local stable manifold $W^s_{loc}(x)$ and a local unstable
manifold $W^u_{loc}(x)$, and two Cantor sets $C^s \subset
W^s_{loc}(x)$ and $C^u \subset W^u_{loc}(x)$, such that, for any
$y_s \in C^s$ and $y_u \in C^u$, then $y_s$ has a local unstable
manifold $W^u_{loc}(y_s)$ and $y_u$ has a local stable manifold
$W^s_{loc}(y_u)$. Moreover, these two local manifolds intersect at
exactly one point, and this point belongs to $R$.

An $s$-subrectangle of $R$ is a set $\left(\bigcup_{y\in C}
W^s_{loc}(y) \right) \cap R$, where $C$ is a subset of $C_u$. A
$u$-subrectangle is defined in the same way.

\cite{markarian:slow} and \cite{chernov:slow} have proved that
$T:X \to X$ satisfies Chernov's axioms of \cite{chernov:decay}.
This implies that it admits a hyperbolic Young tower in the
following sense: there exist a rectangle $R$ of positive measure,
a partition $R=\bigcup R_i$ (modulo $0$) by $s$-subrectangles, and
return times $r_i\in \N$ such that $T^{r_i}$ is a
homeomorphism on $R_i$, and $T^{r_i}(R_i)$ is a
$u$-subrectangle of $R$. Moreover, the tails of the tower are
exponentially small: there exist $\rho<1$ and $C>0$ such that
  \begin{equation*}
  \forall n\in \N, \mu\left( \bigcup_{r_i>n} R_i \right) \leq C
  \rho^n.
  \end{equation*}

We can then define an abstract space $\bar \Delta$ as the disjoint
union of the sets $T^k (R_i)$ for $i\in \N$ and $k<r_i$. It is
endowed with a natural projection $\pi_X : \bar \Delta \to X$ and
a dynamics $\bar U:\bar \Delta \to \bar \Delta$ such that $\pi_X
\circ \bar U= T \circ \pi_X$.

It is also possible to construct on $\bar \Delta$ a probability
measure $\mu_{\bar \Delta}$ which is invariant under $\bar U$ and
such that $(\pi_X)_*( \mu_{\bar \Delta}) = \mu$. Note however that
$\pi_X$ is in general strongly not injective, so that $\mu_{\bar
\Delta}$ can not be defined as the pullback of $\mu$. Rather, one
constructs an invariant measure for $\bar U$, and one proves that
its projection, being absolutely continuous with respect to $\mu$
and $T$-invariant, is necessarily $\mu$.

It is then useful to go from this abstract hyperbolic dynamics to
an abstract expanding dynamics. To do so, one identifies the
points of $\bar \Delta$ which are on the same stable leaf in some
rectangle. This defines a space $\Delta$, together with a
projection $\pi_\Delta: \bar \Delta \to \Delta$. Since the map
$\bar U$ sends stable leaves to stable leaves, it gives rise to a
dynamics $U:\Delta \to \Delta$ on the quotient. The measure
$\mu_\Delta:= (\pi_\Delta)_*( \mu_{\bar\Delta})$ is invariant
under $U$. Then $(\Delta, U,\mu_\Delta)$ is an \emph{expanding
Young tower}, in the sense of Section \ref{subs:towers}.

\subsection{Coboundary results}
\label{sec:coboundary}

Let $f_0:X_0 \to \R$ be a H\"older function satisfying $(P1)$, for
which we want to prove a limit theorem. Since it is easier to work
in an expanding and well understood setting, we will first prove
results in $\Delta$, and then go back from $\Delta$ to $X_0$.

For $x\in X$, let first $f(x)=\sum_{k=0}^{\phi_+(x)-1}f_0(T_0^k
x)$. This function is not bounded any more. However, if two points
$x$ and $y$ are on a local stable manifold which is not cut by a
discontinuity of $T_0$ during the next $n$ iterates of $T_0$, and
with $\phi_+(x)=\phi_+(y)=n$, then
  \begin{equation}
  \label{contracts}
  |f(x)-f(y)| \leq C n d(x,y)^\alpha
  \end{equation}
for some $\alpha>0$. In the same way, if $x$ and $y$ are on a
local unstable manifold which is not cut during the next $n$
iterates of $T_0$, and $\phi_+(x)=\phi_+(y)=n$, then $|f(x)-f(y)|
\leq C n d(T x, T y)^\alpha$. Moreover, the property $(P1)$ implies
that $f$ is bounded on the set of points sliding along the
semicircles.

The measure of points bouncing $n$ times along the segments is
$\sim \frac{\ell^2}{\pi n^3}$. Since $I$ (given by \eqref{donne_I})
is nonzero, the function $f$ is equivalent to $nI$ on this set,
and we obtain
  \begin{equation}
  \label{domaine_attraction}
  \mu \{ x\tq |f(x)| \geq n\} \sim \sum_{n/|I|}^\infty \frac{\ell^2}{\pi
  k^3} \sim \frac{I^2 l^2}{2\pi n^2}.
  \end{equation}
Hence, the distribution of $f$ is in the nonstandard domain of
attraction of the Gaussian law (see
Paragraph~\ref{par_attraction}).

Define a function $\bar f$ on $\bar \Delta$ by $\bar f = f\circ
\pi_X$. It would be easy to go finally from $\bar \Delta$ to
$\Delta$ if $\bar f$ were constant along the local stable leaves
in $\bar \Delta$ (which would mean that $\bar f$ would induce a
function on the quotient $\Delta$). This is in general not the
case, but we will prove that $\bar f$ is cohomologous to such a
function, using the usual cohomology trick.

For every rectangle in $\bar \Delta$, choose a definite unstable
leaf. Define a projection $\pi : \bar \Delta \to \bar \Delta$ by
sliding along stable manifolds to this specific unstable manifold.
We define a function $\bar u(x)=\sum_{k=0}^\infty \bigl[ \bar
f(\bar U^k x) - \bar f(\bar U^k \pi x)\bigr]$. Note that, despite
of the fact that $T$ contracts stable manifolds uniformly, the
function ${\bar u}(x)$ may not seem well-defined at first sight,
as ${\bar f}$ -- and consequently, its H\"older constant -- is
unbounded. Nevertheless, whenever $\bar f$ is large, $T$ contracts
stable manifolds strongly, and the H\"older constant can be regained
by going down the tower. This is the essence of the following
lemma.

\begin{lem}
\label{u_bornee}
The function $\bar u$ is well defined and bounded on $\bar
\Delta$.
\end{lem}
\begin{proof}
In this proof the positive constants $C$ do depend on the H\"older
exponent $\alpha$, but this has no significance. Let $K\in \N$ be
such that $\alpha K>1$. Consider first $x$ which is at height
$\geq K$ in the tower. Let $y=\pi x$. Let $x'=\bar U^{-K}x$ and
$y'=\bar U^{-K} y$. We will prove that
  \begin{equation}
  \label{a_prouver}
  \forall k\in \N, |\bar f(\bar U^k x)-\bar f(\bar U^k y)|
  \leq C d(\pi_X \bar U^k x', \pi_X \bar U^k y')^\alpha.
  \end{equation}
Namely, if $\phi_+( \pi_X \bar U^k x )=n$, then
  \begin{equation}
  \label{estime_fbar}
  |\bar f(\bar U^k x)-\bar f(\bar U^k y)|
  \leq C n d(\pi_X \bar U^k x, \pi_X \bar U^k y)^\alpha= C
  n d( \pi_X \bar U^{k+K}x', \pi_X \bar U^{k+K}y')^\alpha
  \end{equation}
by \eqref{contracts}. If $n=\phi_+(\pi_X \bar U^k x)$ is bounded,
the conclusion is trivial. If $n$ is large and $\pi_X \bar U^k x$
is a sliding point, the conclusion is also trivial by $(P1)$.

Hence, assume that $n$ is large and that $\pi_X \bar U^k x$ is a
bouncing point. Proposition \ref{prop_contracts} implies that, for
$0 \leq i< K$, $\phi_+(\pi_X \bar U^{k+i} x') \geq n/4^{K-i} \geq
n/4^K$. Once again by Proposition  \ref{prop_contracts}, we get
  \begin{equation*}
  d(\pi_X \bar U^{k+i+1}x', \pi_X \bar U^{k+i+1}y') \leq \frac{C}{n} d(\pi_X
  \bar U^{k+i}x', \pi_X \bar U^{k+i}y').
  \end{equation*}
Hence,
  \begin{equation*}
  d( \pi_X \bar U^{k+K}x', \pi_X \bar U^{k+K}y')
  \leq \frac{C}{n^K}d(\pi_X \bar U^k x', \pi_X \bar U^k y').
  \end{equation*}
Together with \eqref{estime_fbar} and the inequality $K\alpha>1$,
this implies \eqref{a_prouver}.

Since $\pi_X x'$ and $\pi_X y'$ are on a local stable manifold,
$d(\pi_X \bar U^k x', \pi_X \bar U^k y')$ goes exponentially fast
to zero. Hence, the series $\sum |\bar f(\bar U^k x) -\bar f(\bar
U^k y)|$ is summable, and $\bar u(x)$ is well defined.

Suppose now that $x$ is at height $<K$. Let $y=\pi x$. Applying
the previous argument to $x'=x$ and $y'=y$, we get that
$\sum_{k=K}^\infty |\bar f(\bar U^k x)-\bar f(\bar U^k y)|$ is
bounded. Moreover, during the first $K$ iterates, $x$ and $y$
remain at a bounded height in the tower, which implies that $\bar
f(\bar U^k x)$ and $\bar f(\bar U^k y)$ remain uniformly bounded.
This concludes the proof.
\end{proof}

Let $\bar g(x)=\bar f(x)-\bar u(x)+\bar u(\bar U x)$. Then
  \begin{equation*}
  \bar g(x)=\bar f(\pi x) + \sum_{k=0}^\infty \Bigl[\bar f(\bar U^k(\bar U
  \pi x)) -\bar f( \bar U^k (\pi \bar U \pi x))\Bigr].
  \end{equation*}
Hence, $\bar g(x)$ depends only on $\pi x$, i.e., $\bar g$ is
constant along the stable manifolds in the rectangles.
Consequently, there exists a function $g:\Delta\to \R$ such that
$\bar g=g\circ \pi_\Delta$.

It will be important that $g$ is regular enough on $\Delta$, to
use functional analytic techniques. For $x_1,x_2\in \Delta$, let
$s(x_1,x_2)$ be their separation time, i.e., the number of returns
to the basis before $x_1$ and $x_2$ get into different elements of
the partition. To obtain the following lemma, we will use several
times the same argument as in Lemma~\ref{u_bornee}, but sometimes
along unstable manifolds instead of stable ones.

\begin{lem}
\label{lem:g_Holder}
There exist $C>0$ and $\tau<1$ such that, for every $x_1,x_2$ in
the same element of partition of $\Delta$,
  \begin{equation*}
  | g(x_1)-g(x_2)| \leq C \tau^{s(x_1,x_2)}.
  \end{equation*}
\end{lem}
\begin{proof}
Let us first prove that, if $x_1,x_2$ belong to the same unstable
leaf in a rectangle of $\bar \Delta$, then
  \begin{equation}
  \label{g_unif_bounded}
  \bar g(x_1)-\bar g(x_2) \text{ is uniformly bounded.}
  \end{equation}
The same argument as in the proof of Lemma~\ref{u_bornee} shows
that $\sum_{k=0}^\infty \bigl[\bar f(\bar U^k(\bar U
  \pi x)) -\bar f( \bar U^k (\pi \bar U \pi x))\bigr]$ is bounded. Hence, it
is sufficient to prove that $\bar f(\pi x_1)-\bar f(\pi x_2)$ is
bounded. Let $K$ be as in the proof of Lemma~\ref{u_bornee}. If
$x_1$ (and $x_2$) return to the basis of $\bar \Delta$ before time
$K$, then $\phi_+(x_1)=\phi_+(x_2)$ is bounded, which implies that
$\bar f(\pi x_1)$ and $\bar f(\pi x_2)$ are bounded. If $x_1$ (and
$x_2$) are sliding points, then the conclusion is also a
consequence of $(P1)$. Otherwise, $x_1$ and $x_2$ are bouncing
points. We show as in the proof of Lemma~\ref{u_bornee} (but along
the \textit{unstable} leaf containing $\pi x_1$ and $\pi x_2$)
that $|\bar f(\pi x_1) -\bar f(\pi x_2)| \leq C d(\pi_X \bar U^K
\pi x_1, \pi_X \bar U^K \pi x_2)^\alpha$. Since this quantity is
uniformly bounded, this concludes the proof of
\eqref{g_unif_bounded}.

Take $x_1,x_2 \in \bar \Delta$ on the same unstable leaf, and let
$s=s(\pi_\Delta (x_1), \pi_\Delta(x_2))$. We will prove that
  \begin{equation}
  \label{g_Lipschitz}
  |\bar g(x_1) -\bar g(x_2)| \leq C \lambda^{\alpha s/2}
  \end{equation}
for some $C>0$, where $\lambda<1$ is larger than the contraction
coefficient of $T$ along stable manifolds, and the contraction
coefficient of $T^{-1}$ along unstable manifolds.

By \eqref{g_unif_bounded}, this is trivial if $s<2K$. Hence, we
can assume $s\geq 2K$. Let $N=\left\lfloor
\frac{s}{2}\right\rfloor\geq K$, then
  \begin{equation}
  \label{estime_g_4termes}
  \begin{split}
  \bar g(x_1)-\bar g(x_2) & =
 \bar f(\pi x_1) -\bar f(\pi x_2)
  \\&
  +\sum_{k=0}^{N-1} \Bigl[ \bar f(\bar U^k(\bar U \pi x_1)) -
  \bar f(\bar U^k(\bar U \pi x_2)) \Bigr]
  +\sum_{k=0}^{N-1} \Bigl[ \bar f(\bar U^k(\pi \bar U \pi x_2)) -
  \bar f(\bar U^k(\pi \bar U \pi x_1)) \Bigr]
  \\ &
  + \sum_{k=N}^\infty \Bigl[ \bar f(\bar U^k (\bar U \pi x_1)) - \bar
  f(\bar U^k (\pi \bar U \pi x_1)) \Bigr]
  + \sum_{k=N}^\infty \Bigl[ \bar f(\bar U^k (\pi \bar U \pi x_2)) - \bar
  f(\bar U^k (\bar U \pi x_2)) \Bigr].
  \end{split}
  \end{equation}
Since $N+K \leq s$, we have for any $k< N$
  \begin{align*}
  \left|\bar f(\bar U^k(\bar U \pi x_1)) -
  \bar f(\bar U^k(\bar U \pi x_2)) \right|&
  \leq C d(\pi_X \bar U^{k+K} (\bar U\pi x_1),
  \pi_X \bar U^{k+K} (\bar U \pi x_2))^\alpha
  \\&
  \leq C \lambda^{\alpha(s-(k+K+1))}
  d(\pi_X \bar U^s \pi x_1, \pi_X \bar U^s \pi x_2)^\alpha
  \leq C \lambda^{\alpha(s-k)}.
  \end{align*}
Summing over $k$, we obtain
  \begin{equation*}
  \left| \sum_{k=0}^{N-1} \Bigl[ \bar f(\bar U^k(\bar U \pi x_1)) -
  \bar f(\bar U^k(\bar U \pi x_2)) \Bigr] \right|
  \leq C \lambda^{\alpha(s-N)} \leq C \lambda^{\alpha s/2}.
  \end{equation*}
The other term on the second line of \eqref{estime_g_4termes} can
be estimated in the same way, as well as the term on the first
line of \eqref{estime_g_4termes}.

Since $N \geq K$, we also have for any $k\geq N$
  \begin{align*}
  \left| \bar f(\bar U^k (\bar U \pi x_1)) - \bar
  f(\bar U^k (\pi \bar U \pi x_1)) \right|
  &
  \leq C d(\pi_X \bar
  U^{k-K}(\bar U \pi x_1), \pi_X \bar U^{k-K} (\pi \bar U \pi
  x_1))^\alpha
  \\&
  \leq C \lambda^{\alpha(k-K)} d(\pi_X (\bar U \pi x_1), \pi_X (\pi
  \bar U\pi x_1))^\alpha
  \leq C \lambda^{\alpha k}.
  \end{align*}
Summing over $k$, we obtain
  \begin{equation*}
  \left|\sum_{k=N}^\infty \Bigl[ \bar f(\bar U^k (\bar U \pi x_1)) - \bar
  f(\bar U^k (\pi \bar U \pi x_1)) \Bigr] \right|
  \leq C \lambda^{\alpha N} \leq C \lambda^{\alpha s/2}.
  \end{equation*}
The other term on the third line of \eqref{estime_g_4termes} is
handled in the same way.
\end{proof}

\section{Limit theorems in Young towers}

\subsection{A result by Aaronson and Denker}

\label{par_attraction}

A function $f:\R_+^* \to \R_+^*$ is \emph{slowly varying} if, for all
$\lambda>0$, $f(\lambda x)/f(x)$ tends to $1$ when $x\to\infty$.

By classical probabilistic results,
a real random variable $Z$ is in the nonstandard domain of
attraction of the Gaussian distribution $\boN(0,1)$ if and only if
it satisfies one of the following equivalent conditions:
\begin{itemize}
\item The function $L(x):=E(Z^2 1_{|Z|\leq x})$ is unbounded
and slowly varying.
\item $P(|Z|>x) \sim x^{-2} l(x)$ for some function $l$ such that
$\tilde{L}(x):= 2\int_1^x \frac{l(u)}{u}\dd u$ is unbounded and
slowly varying.
\end{itemize}

\begin{rmq}
In this case, $\tilde{L}(x) \sim L(x)$ when $x \to \infty$, and
$l(x)=o(L(x))$. It is possible, however, that $l$ is not slowly
varying and that these conditions hold anyway.
\end{rmq}

Such a random variable belongs to $L^p$ for all $1\leq p<2$, but
not to $L^2$. We will say that $l$ and $L$ are the \emph{tail
functions} of $Z$. They are defined up to asymptotic equivalence.
Choose a sequence $B_n \to \infty$ such that $\frac{n}{B_n^2}
L(B_n) \to 1$. Then, if $Z_0,Z_1,\dots$ is a sequence of
independent random variables distributed as $Z$, then
  \begin{equation*}
  \frac{Z_0+\dots+Z_{n-1} - nE(Z)}{B_n} \to \boN(0,1).
  \end{equation*}

More generally, if $\frac{n}{B_n^2}L(B_n) \to C>0$, then the
previous sequence converges to $\boN(0,C)$.

In \cite{aaronson_denker:central}, Aaronson and Denker have proved
the same kind of limit theorem when the sequence $Z_0,Z_1,\dots$
is not independent. More precisely, consider $U$ a mixing
Gibbs-Markov map (as defined in \cite{aaronson:book}) on a space
$\Delta$, preserving a probability measure $\mu_\Delta$, and let
$g:\Delta \to \R$ be a function which is locally H\"older and whose
distribution with respect to $\mu_\Delta$ is in the nonstandard
domain of attraction of $\boN(0,1)$ as above. Then they prove that
  \begin{equation*}
  \frac{g+g\circ U+\dots +g\circ U^{n-1} -n \int g}{B_n} \to \boN(0,1)
  \end{equation*}
as above.

The proof goes as follows: let $\hat{U}$ be the transfer operator
associated to $U$, and $\hat{U}_t$ its perturbation given by
$\hat{U}_t
 u=\hat{U}( e^{itg} u)$. These operators satisfy a Lasota-Yorke
inequality on the space of H\"older functions, and $\bigl\|
\hat{U}_t-\hat{U} \bigr\|=O(t)$. Hence, the eigenvalue $\lambda_t$
of $\hat{U}_t$ close to $1$ satisfies $|\lambda_t-1|=O(t)$, and
the corresponding eigenfunction $w_t$ (normalized so that $\int
w_t=1$) is such that $\norm{w_t-1}=O(t)$.

Then they prove the following abstract lemma (in which there is no
dynamics, this lemma depends only on the distribution of $g$):
\begin{lem}
\label{lem_aaronson_denker}
For any bounded function $w$ on $\Delta$,
  \begin{equation*}
  \int (e^{itg}-1-itg) w  = -\frac{t^2}{2}
  \int 1_{|g| \leq 1/|t|} g^2 w  +
\norm{w}_\infty o(t^2L(1/|t|)).
  \end{equation*}
Here, the $o(t^2 L(1/|t|))$ is uniform in $w$.
\end{lem}

Applying this lemma to $w_t$, one gets
  \begin{equation*}
  \lambda_t-1-it\int gw_t= \int (e^{itg}-1-itg) w_t
  = -\frac{t^2}{2} \int 1_{|g| \leq 1/|t|} g^2 w_t  +  o(t^2L(1/|t|))
  \end{equation*}
(where we have used the fact that $w_t$ is bounded). Since
$\norm{w_t-1}_\infty =o(1)$,
  \begin{equation*}
  \frac{t^2}{2} \int 1_{|g| \leq 1/|t|} g^2 w_t
  =\frac{t^2}{2} L(1/|t|) (1+o(1)).
  \end{equation*}
Hence,
  \begin{equation}
  \label{ad_final}
  \lambda_t=1+it \int gw_t - \frac{t^2}{2} L(1/|t|) (1+o(1)).
  \end{equation}
Finally, $\int gw_t=\int g +O(t)$ since
$\norm{w_t-1}_\infty=O(t)$. So we get
  \begin{equation*}
  \lambda_t=1+it \int g - \frac{t^2}{2} L(1/|t|) (1+o(1)).
  \end{equation*}
This expansion is sufficient to get the required limit theorem.

\subsection{The result in Young towers}

\label{subs:towers}

Let $(\Delta,\mu_\Delta)$ be a probability space and $U:\Delta \to
\Delta$ a probability preserving map. We say that $(\Delta,U)$ is
an \emph{expanding Young tower}
 (\cite{lsyoung:recurrence})
if there exist integers $r_p\in \N^*$ and a partition
$\{\Delta_{k,p}\}_{p\in \N, k\in \{0,\dots, r_p-1\}}$ of $\Delta$
such that
\begin{enumerate}
\item
For all $p$ and $k<r_p-1$, $T$ is a measurable isomorphism between
$\Delta_{k,p}$ and $\Delta_{k+1,p}$, preserving $\mu_\Delta$.
\item For all $p$,
$T$ is a measurable isomorphism between $\Delta_{r_p-1,p}$ and
$\Delta_0:= \bigcup_m \Delta_{0,m}$.
\item Let $U_0$ be the first return map induced by $U$ on $\Delta_0$.
For $x,y\in \Delta_0$, define their separation time $s(x,y)=\inf\{
n\in \N \tq U_0^n(x) \text{ and }U_0^n(y)\text{ are not in the
same }\Delta_{0,p}\}$. We extend this separation time to the whole
tower in the following way: if $x,y$ are not in the same element
of partition, set $s(x,y)=0$. Otherwise, $x,y\in \Delta_{k,p}$.
Let $x',y'\in \Delta_{0,p}$ be such that $x=U^k x'$ and $y=U^k
y'$, and set $s(x,y)=s(x',y')$.

For $x\in \Delta$, let $J(x)$ be the inverse of the jacobian of
$U$ at $x$. We assume that there exist $\beta<1$ and $C>0$ such
that, for all $x,y$ in the same element of partition,
  \begin{equation}
  \label{dist_bornee}
  \left| 1-\frac{J(x)}{J(y)} \right| \leq C \beta^{s(Ux,Uy)}
  \end{equation}
\end{enumerate}

\begin{rmq}
Note that the definition of separation time in
\cite{lsyoung:annals} is in terms of the number of all iterations
of $U$, while we follow the convention of
\cite{lsyoung:recurrence} when we define separation in terms of
returns to the basis. Hence, our setting is more general than that
of \cite{lsyoung:annals}, but it will make the proof of the
spectral gap more complicated.
\end{rmq}

Let $\Delta_n=\bigcup \Delta_{n,p}$. This is the set of points at
height $n$ in the tower.
We will say that $(U,\Delta)$ is an \emph{expanding Young tower
with exponentially small tail} if there exists $\rho<1$ such that
$\mu_{\Delta}( \Delta_n) = O(\rho^n)$.

Let $J^{(n)}$ be the inverse of the jacobian of $U^n$. It is
standard that \eqref{dist_bornee} implies that the distortion of
the iterates of $U$ is uniformly bounded, in the following sense:
there exists $C>0$ such that, for all points $x,y$ such that $T^k$ and 
$T^k y$ remain in the same elements of the partition for $0\leq k<n$,
  \begin{equation}
  \label{distn_bornee}
  \left| 1-\frac{ J^{(n)}(x)}{J^{(n)}(y)} \right| \leq C \beta^{s(U^n
  x,U^n y)}.
  \end{equation}

A function $g:\Delta \to \R$ is \emph{locally H\"older} if there
exist $C>0$ and $\tau<1$ such that $|g(x)-g(y)|\leq C
\tau^{s(x,y)}$ for all $x,y$ in the same element of the partition.
This is exactly the type of functions that arise from the stadium
billiard, cf.\ Lemma~\ref{lem:g_Holder}. Note that $g$ can very
well be unbounded. Without loss of generality, we can assume
$\tau\geq \beta$.

Let $\haut(x)$ be the height of the point $x$, i.e., $\omega(x)=n$ if
$x\in \Delta_n$. Let $\pi_0 :
\Delta\to \Delta_0$ be the projection to the basis, and define a
function $G$ on $\Delta$ by $G(x)=\sum_{k=0}^{\haut(x)-1} g(U^k
\pi_0 x)$. In this setting, we get the following extension of the
theorem proved by Aaronson and Denker:
\begin{thm}
\label{thm_ad_young}
Let $U:\Delta\to \Delta$ be an expanding Young tower with
exponentially small tail, and let $g:\Delta \to \R$ be locally H\"older
continuous. Assume that the distribution of $g$ is in the
nonstandard domain of attraction of $\boN(0,1)$, with tail
functions $l$ and  $L$. Assume moreover that $l$ and $L$ are
slowly varying, and $l(x \ln x)/l(x) \to 1$, $L(x\ln x)/L(x) \to
1$ when $x\to \infty$. Finally, assume that there exists a real
number $a\not=-1/2$ such that
  \begin{equation}
  \label{condition_integrale}
  \int g (e^{itG}-1) = (a+o(1))it L(1/|t|)) \text{ when }t\to 0.
  \end{equation}
Write $L_1(x)=(2a+1)L(x)$, and choose a sequence $B_n \to \infty$
such that $\frac{n}{B_n^2}L_1(B_n) \to 1$. Then
  \begin{equation*}
  \frac{S_n g - n \int g}{B_n} \to \boN(0,1).
  \end{equation*}
\end{thm}

The additional assumption on $l$ and $L$ is satisfied in most
natural cases (for example when $l=1$ and $L=\ln$, which will be
the case for the stadium billiard).

When $a=0$, we get the same asymptotics as in Aaronson-Denker's
Theorem. However, when $a\not=0$, then there is an additional
effect due to the presence of the tower.

The proof will follow the same lines as in Aaronson-Denker's
proof: it is possible to construct a good space on which the
transfer operator $\hat{U}$ has a spectral gap. The perturbed
operator $\hat{U}_t$ also has a spectral gap, which gives an
eigenvalue $\lambda_t$ and an eigenfunction $w_t$. The main
problem is that $\bigl\| \hat{U}_t - \hat{U}\bigr\|$ can not be
$O(t)$ in general: it is easy to construct examples where
$t=o\Bigl(\bigl\| (\hat{U}_t-\hat{U})1 \bigr\|_{L^2}\Bigr)$,
whence $t=o\Bigl( \bigl\| \hat{U}_t - \hat{U}\bigr\|\Bigr)$ as
soon as the good space is contained in $L^2$ and contains the
function $1$.

Using abstract arguments by Keller and Liverani, we can
nevertheless prove that $|\lambda_t -1|=O(|t|^{1/10})$ and
$\norm{w_t-1}_{L^1}=O(|t|^{1/10})$. These information are
(essentially) sufficient to apply Aaronson and Denker's argument
and get $\lambda_t=1+it\int gw_t -\frac{t^2}{2}\int 1_{|g|\leq
1/|t|} g^2 w_t+ o(t^2 L(1/|t|))$ as in \eqref{ad_final}. The main
difficulty is then to make the function $w_t$ disappear in this
expression, to get something more tractable. We will namely show
that $\int 1_{|g|\leq 1/|t|} g^2 w_t \sim L(1/|t|)$ and $\int gw_t
\sim \int g e^{itG}$, which will conclude the proof.

To do this, we need to know that $w_t-1=O(t)$ in some sense. To
prove such an estimate, we use a roundabout technical argument
relying on the fact that the induced map on the basis of the tower
is uniformly expanding, to prove that $\norm{
1_{\Delta_0}(w_t-1)}_\infty=O(t)$, and then we propagate this
information up in the tower, using the information we have already
proved on $\lambda_t$. This propagation requires the Birkhoff sums
of $g$ to be small enough. To ensure this on a set of large
measure, we use the information on the tails of $g$. This is the
only point where the additional information on $l$ and $L$ is
used.

\subsection{Proof of Theorem~\ref{thm_ad_young}}

We will first prove Theorem~\ref{thm_ad_young} assuming that $\int
g=0$. In Paragraph~\ref{etend_cas_general}, we will show that this
implies the theorem in full generality. Hence, \emph{until the end
of Paragraph~\ref{cas_particulier}, we will assume that $\int
g=0$}.

\subsubsection{Construction of the functional spaces and the transfer
operators}

Since the tails of the tower are exponentially small by
assumption, there exists $\rho<1$ such that $\mu_\Delta(\Delta_n)
\leq C \rho^n$. Denote the return time to the basis from itself by
$\phi$. Take $\epsilon>0$ such that $e^{6\epsilon}\rho<1$.

For $u:\Delta\to \C$, write
 \begin{equation*}
 \norm{u}_{m}=\inf\{C \tq \forall n\in \N,
  \text{for almost every }x\in \Delta_n, |u(x)|\leq C e^{\epsilon n}\}
 \end{equation*}
and
 \begin{multline*}
 \norm{u}_l = \inf\{C \tq \text{for almost every }x,y\text{ in the same
element of the partition at height }n,\\ |u(x)-u(y)|\leq
Ce^{\epsilon n} \tau^{s(x,y)}\}.
  \end{multline*}

Denote by $\boH$ the space of measurable functions $u$ on $\Delta$
for which $\norm{u}:=\norm{u}_m + \norm{u}_l <+\infty$. It is a
Banach space included in $L^1$ (and even in $L^6$ because of the
condition $e^{6\epsilon}\rho<1$). This inclusion is compact.

The following proposition is similar to a result of Young:
\begin{prop}
\label{ls_hatU}
There exist $C>0$ and $\theta<1$ such that, for any $u\in \boH$,
for any $n\in \N$,
  \begin{equation*}
  \bigl\| \hat{U}^n u \bigr\| \leq C\theta^n \norm{u} + C \norm{u}_{L^1}.
  \end{equation*}
\end{prop}
Note that our definition of separation time is not the same as in
\cite{lsyoung:annals}, and that Young uses the fact that the
return to the basis only occur after a large time $N$. This gives
her a strong expansion, sufficient to get rid of constants easily.
This is not true in our setting. Hence, the proof of the
proposition will be more involved than Young's.

\begin{proof}
Take $x\in \Delta_0$. Then $\hat{U}^n u(x)=\sum
J^{(n)}(x_p)u(x_p)$, where $\{x_p\}=U^{-n}(x)$. Let $A_p$
containing $x_p$ be such that $U^n : A_p \to \Delta_0$ is an
isomorphism. Then $J^{(n)}(x_p) \leq C \mu_\Delta(A_p)$ since the
distortion is bounded, by \eqref{distn_bornee}. 
Let $\omega_p$ be the height of the set
$A_p$ and $r_p$ the number of returns of $A_p$ to the basis before
time $n$.

For $y \in A_p$, $s(x_p,y) \geq r_p$, whence
 \begin{equation*}
 |u(x_p)-u(y)| \leq \tau^{r_p} e^{\epsilon \omega_p} \norm{u}_l.
  \end{equation*}
Hence,
  \begin{equation}
  \label{borne_uxp}
  |u(x_p)| \leq \tau^{r_p} e^{\epsilon \omega_p} \norm{u}_l +
  \frac{1}{\mu_\Delta(A_p)}\int_{A_p} |u|.
  \end{equation}
We get
  \begin{equation}
  \label{eq_ly_normem}
  |\hat{U}^n u(x)| \leq C \sum \mu_\Delta(A_p) \tau^{r_p} e^{\epsilon \omega_p}
  \norm{u}_l + C \int |u|.
  \end{equation}
Let $\omega: \Delta \to \N$ be the function ``height'', and let
$\Psi_n(x)$ be the number of returns of $x$ to the basis between
time $1$ and $n$. Then \eqref{eq_ly_normem} implies that
  \begin{equation}
  \label{eq_ly_normem2}
  |\hat{U}^n u(x)| \leq C \norm{u}_l \int_{U^{-n}\Delta_0} \tau^{\Psi_n}
  e^{\epsilon \omega} + C \norm{u}_{L^1}.
  \end{equation}

We will use the following technical lemma, which will be proved in
the appendix.
\begin{lem}
\label{renouvell_exp}
There exist $C>0$ and $\theta<1$ such that, for any $n\in \N$,
  \begin{equation}
    \int_{U^{-n}\Delta_0} \tau^{\Psi_n}
  e^{\epsilon \omega} \leq C \theta^n.
  \end{equation}
\end{lem}
Increasing $\theta$ if necessary, we can assume that
$e^{-\epsilon} \leq \theta$.

This lemma, together with \eqref{eq_ly_normem2}, implies that, for
any $x\in \Delta_0$,
  \begin{equation}
  |\hat{U}^n u(x)|\leq C \theta^n \norm{u}_l +C \norm{u}_{L^1}.
  \end{equation}

Consider now $x\in \Delta$ such that $\haut(x)=k<n$. Let $x'$ be
its projection in the basis. Then $\hat{U}^n
u(x)=\hat{U}^{n-k}u(x')$, whence
  \begin{equation}
  e^{-\epsilon k} |\hat{U}^n u(x)|
  =e^{-\epsilon k} |\hat{U}^{n-k}u(x')|
  \leq e^{-\epsilon k} C \theta^{n-k} \norm{u}_l + Ce^{-\epsilon k}
\norm{u}_{L^1}
  \leq C \theta^n  \norm{u}_l + C \norm{u}_{L^1}.
  \end{equation}
Assume finally that $\haut(x)=k \geq n$. Let $x'=U^{-n}(x)$. Then
  \begin{equation}
  e^{-\epsilon k} |\hat{U}^n u(x)|=e^{-\epsilon n} e^{-\epsilon(k-n)}|u(x')|
  \leq e^{-\epsilon n} \norm{u}_m.
  \end{equation}

These equations prove that
  \begin{equation*}
  \bigl\| \hat{U}^n u\bigr\|_m \leq C \theta^n \norm{u}+C \norm{u}_{L^1}.
  \end{equation*}

We still have to handle the H\"older norm. Consider two points $x,y$
in the same element of partition of the basis $\Delta_0$. Let
$x_p$ and $y_p$ be their preimages, in sets $A_p$ as above. Then
  \begin{align*}
  |\hat{U}^n u(x) - \hat{U}^nu(y)|&
  \leq \sum |J^{(n)}(x_p)u(x_p)-J^{(n)}(y_p)u(y_p)|
  \\&
  \leq \sum |J^{(n)}(x_p)| |u(x_p)-u(y_p)| + \sum |J^{(n)}(x_p)|
  \left| 1-\frac{J^{(n)}(y_p)}{J^{(n)}(x_p)}\right| |u(y_p)|.
  \end{align*}
In the first sum, $|J^{(n)}(x_p)| \leq C \mu_\Delta(A_p)$ and
$|u(x_p)-u(y_p)| \leq \tau^{s(x,y)+r_p}e^{\epsilon \omega_p}
\norm{u}_l$. Hence,
  \begin{equation*}
  \sum |J^{(n)}(x_p)| |u(x_p)-u(y_p)| \leq C \tau^{s(x,y)}\norm{u}_l
  \int_{U^{-n}\Delta_0} \tau^{\Psi_n} e^{\epsilon \omega}
  \leq C \theta^n \tau^{s(x,y)}\norm{u}_l
  \end{equation*}
by Lemma~\ref{renouvell_exp}.

In the second sum, $|J^{(n)}(x_p)| \leq C \mu_\Delta(A_p)$ and
$\left| 1-\frac{J^{(n)}(y_p)}{J^{(n)}(x_p)}\right| \leq C
\tau^{s(x,y)}$ by \eqref{distn_bornee}. Moreover, $|u(y_p)|$ is
bounded by \eqref{borne_uxp}. Using these inequalities, we get
  \begin{align*}
  \sum |J^{(n)}(x_p)|
  \left| 1-\frac{J^{(n)}(y_p)}{J^{(n)}(x_p)}\right| |u(y_p)|&
  \leq \sum C \mu_\Delta(A_p) \tau^{s(x,y)} \left[ \tau^{r_p} e^{\epsilon \omega_p}
  \norm{u}_l +
  \frac{1}{\mu_\Delta(A_p)}\int_{A_p} |u| \right]
  \\&
  \leq C \tau^{s(x,y)} \norm{u}_l \int_{U^{-n}\Delta_0}
  \tau^{\Psi_n}e^{\epsilon \omega}+ C \tau^{s(x,y)}\int |u|
  \\&
  \leq  C \theta^n \tau^{s(x,y)}\norm{u}_l + C \tau^{s(x,y)} \norm{u}_{L^1}
  \end{align*}
by Lemma~\ref{renouvell_exp}.

To sum up, we have proved that, when $x$ and $y$ belong to the
same partition element of the basis,
  \begin{equation*}
  \frac{|\hat{U}^n u(x) - \hat{U}^n u(y)|}{\tau^{s(x,y)}} \leq C
  \theta^n \norm{u}_l
  +C \norm{u}_{L^1}.
  \end{equation*}

Let now $x$ and $y$ belong to the same element of the partition,
with $k=\haut(x)<n$. Let $x'$ and $y'$ be their projection in the
basis. Then
  \begin{align*}
  e^{-\epsilon k}\frac{|\hat{U}^nu( x)-\hat{U}^n u(y)|}{\tau^{s(x,y)}}&
  =e^{-\epsilon k}\frac{|\hat{U}^{n-k}u(x') - \hat{U}^{n-k}u(y')|}
  {\tau^{s(x',y')}}
  \leq e^{-\epsilon k} \left[C\theta^{n-k} \norm{u}_l + C \norm{u}_{L^1}\right]
  \\&
  \leq C \theta^n \norm{u}_l+C \norm{u}_{L^1}.
  \end{align*}
Assume finally that $k\geq n$. Let $x'=U^{-n}x$ and $y'=U^{-n}y$.
Then
  \begin{equation*}
  e^{-\epsilon k} \frac{ |\hat{U}^nu( x)-\hat{U}^nu (y)|}{\tau^{s(x,y)}}
  =e^{-\epsilon n} e^{-\epsilon(k-n)} \frac{|u(x')-u(y')|}{\tau^{s(x',y')}}
  \leq \theta^n \norm{u}_l.
  \end{equation*}
Summing up these equations, we get $\bigl\|\hat{U}^n u\bigr\|_l
\leq C \theta^n \norm{u}_l + C \norm{u}_{L^1}$. This concludes the
proof of the proposition.
\end{proof}

Let $g$ be the locally H\"older function for which we want to prove
a limit theorem. It is possible that $g\not\in \boH$, since
$\norm{g}_m$ is not necessarily finite.

Define a perturbed transfer operator, à la Nagaev, by $\hat{U}_t
(u) =\hat{U}(e^{itg}u)$.
\begin{prop}
\label{ly_perturb}
There exist constants $C>0$ and $\theta <1$ such that, for all
$t\in [-1,1]$, for all $u\in \boH$, for all $n\in \N$,
  \begin{equation*}
  \bigl\| \hat{U}_t^n u \bigr\| \leq C \theta^n \norm{u}+C \norm{u}_{L^1}.
  \end{equation*}
\end{prop}
This proposition contains Proposition~\ref{ls_hatU} as a special
case, for $t=0$.
\begin{proof}
Let $x\in \Delta$. Then $\hat{U}_t^n u(x)=\sum_{U^ny=x}
e^{itS_ng(y)}J^{(n)}(y)u(y)$, whence $|\hat{U}_t^n u(x)| \leq
\hat{U}^n |u| (x)$. The bound on $\bigl\| \hat{U}^n |u|\bigr\|_m$
thus implies the required bound on $\bigl\|\hat{U}_t^n
u\bigr\|_m$.

For the H\"older norm, take $x$ and $y$ two points in the same
element of partition. Then, with the notations of the proof of
Proposition~\ref{ls_hatU},
  \begin{align*}
  \left|\hat{U}_t^n u(x) - \hat{U}_t^n u(y)\right|&
  =\left| \sum e^{itS_ng (x_p)}J^{(n)}(x_p) u(x_p)- e^{itS_n g(y_p)}
  J^{(n)}(y_p)u(y_p)
  \right|
  \\&
  \leq \sum | J^{(n)}(x_p)u(x_p)-J^{(n)}(y_p)u(y_p) |
  + \sum |e^{it S_ng(x_p)}-e^{itS_ng(y_p)}| J^{(n)}(x_p)|u(x_p)|.
  \end{align*}
The first sum has already been estimated in the proof of
Proposition~\ref{ls_hatU}. For the second one, $|e^{it S_n
g(x_p)}-e^{it S_n g(y_p)}| \leq nC \tau^{s(x,y)}$. Hence,
Proposition~\ref{ls_hatU} implies that
  \begin{equation*}
  \bigl\| \hat{U}_t^n u \bigr\| \leq C (n+1) \theta^n \norm{u}+ C(n+1)
  \norm{u}_{L^1}.
  \end{equation*}
Choose $N>0$ such that $\bar\theta:= C (N+1) \theta^N <1$.
Iterating the equation $\bigl\| \hat{U}_t^N u \bigr\| \leq
\bar\theta \norm{u}+C \norm{u}_{L^1}$ (and using the fact that
$\bigl\|\hat{U}_t^N u\bigr\|_{L^1} \leq \norm{u}_{L^1}$), we get
  \begin{equation*}
  \bigl\|\hat{U}_t^{nN} u\bigr\| \leq \bar \theta^{n} \norm{u} +
  \frac{C}{1-\bar\theta } \norm{u}_{L^1}.
  \end{equation*}
This implies the conclusion of the proposition, for the constant
$\bar\theta^{1/N}<1$.
\end{proof}

\begin{lem}
\label{continuite_faible}
When $t \to 0$, $\| \hat{U}_t -\hat{U} \|_{\boH \to L^3}
=O(|t|^{1/6})$.
\end{lem}
\begin{proof}
For $u \in \boH$, $(\hat{U}_t-\hat{U})u = \hat{U} ((e^{itg}-1)
u)$. The transfer operator $\hat{U}$ is a contraction in every
$L^p$ space, and in particular in $L^3$. Hence,
  \begin{equation*}
  \bigl\| (\hat{U}_t -\hat{U}) u \bigr\|_{L^3} \leq \norm{ (e^{itg}-1)
  u}_{L^3} \leq \norm{ e^{itg}-1}_{L^6} \norm{u}_{L^6}.
  \end{equation*}
Note that $\norm{u}_{L^6} \leq C \norm{u}$. Hence, $\| \hat{U}_t
-\hat{U} \|_{\boH \to L^3} =O( \norm{ e^{itg}-1}_{L^6} )$. To
estimate this quantity, choose $C>0$ such that, for all $x\in \R$,
$|e^{ix}-1|\leq C|x|^{1/6}$. Then $\int |e^{itg}-1|^6 \leq C \int
|tg| = O(|t|)$. Hence, $\norm{ e^{itg}-1}_{L^6}=O(|t|^{1/6})$.
\end{proof}

\subsubsection{Definition of $\lambda_t$, first estimates}

By Proposition~\ref{ls_hatU} and Hennion's Theorem \cite{hennion},
the operator $\hat{U}: \boH \to \boH$ is quasicompact: outside of
the disk $\{| z|\leq \theta\}$, its spectrum is composed of
discrete eigenvalues of finite multiplicity. In particular, by
ergodicity, $1$ is a simple isolated eigenvalue of $\hat{U}$, with
multiplicity one (and the corresponding eigenfunction is the
constant function $1$).

Lemma~\ref{continuite_faible} is \emph{not} a continuity statement
in $\boH$. However, the operators $\hat{U}$ and $\hat{U}_t$
satisfy a uniform Lasota-Yorke inequality between $\boH$ and
$L^3$, by Proposition~\ref{ly_perturb} (and the fact that
$\norm{u}_{L^1} \leq \norm{u}_{L^3}$). Hence, we can apply the
abstract results of \cite[Corollary 1]{keller_liverani} (following
ideas of \cite{baladi_young}). We get the following:

For small enough $t$, $\hat{U}_t$ has a unique eigenvalue
$\lambda_t$ close to $1$, and it satisfies $|\lambda_t -1| =
O(|t|^{1/10})$. Let $P_t$ be the corresponding spectral
projection. Then $\norm{P_t}_{\boH \to \boH}$ is bounded when
$t\to 0$. Moreover, $\norm{P_t-P_0}_{\boH \to L^3} =
O(|t|^{1/10})$.

\begin{rmq}
Here, $1/10$ could be replaced by any exponent $<1/6$, but any
positive exponent would be sufficient for our purposes.
\end{rmq}

Let $\bar w_t := P_t 1$, and write $w_t = \frac{\bar w_t}{\int
\bar w_t}$. Then $w_t$ is bounded in $\boH$ and
  \begin{equation}
  \label{w_t_proche_1}
  \norm{w_t-1}_{L^3}=O(|t|^{1/10}).
  \end{equation}

\begin{lem}
When $t\to 0$,
  \begin{equation}
  \label{decrit_lambdat_1}
  \lambda_t = 1-\frac{t^2}{2}\int 1_{|g|\leq 1/|t|} g^2 w_t
 + it \int g w_t + o(t^2 L(1/|t|)).
  \end{equation}
\end{lem}
\begin{proof}
By definition, $\hat{U}_t (w_t)=\lambda_t w_t$. Integrating, we
get
  \begin{equation}
  \lambda_t = \int e^{itg} w_t.
  \end{equation}
We want to use Lemma~\ref{lem_aaronson_denker} to estimate this
integral. However, this lemma applies only to bounded functions.
Hence, we will have to modify $w_t$.

Take $x\in \Delta$ with $\haut(x)>0$, and let $x'=U^{-1}(x)$. The
equation $\hat{U}_t w_t =\lambda_t w_t$ implies that $e^{itg(x')}
w_t(x') = \lambda_t w_t(x)$. Hence, $|w_t(x)|= |\lambda_t|^{-1}
|w_t(x')|$. Since $w_t$ is uniformly bounded on the basis of the
tower (since it is bounded in $\boH$), we get
  \begin{equation}
  \label{borne_gt}
  |w_t(x)| \leq C |\lambda_t|^{-\haut(x)}.
  \end{equation}
Define a function $w'_t$ by $w'_t(x)=w_t(x)$ if $\haut(x) \leq
|t|^{-1/10}$ and $w'_t(x)=0$ otherwise. Since $\lambda_t =
1+O(|t|^{1/10})$, \eqref{borne_gt} implies that $w'_t$ is
uniformly bounded.

Lemma~\ref{lem_aaronson_denker} applied to $w'_t$ gives
  \begin{equation}
  \label{sur_w't}
  \int (e^{itg}-1-itg) w'_t = -\frac{t^2}{2} \int 1_{|g| \leq 1/|t|}
  g^2 w'_t + o(t^2 L(1/|t|)).
  \end{equation}
Let us show that this equation is also satisfied by
$w''_t:=w_t-w'_t$. First,
  \begin{equation*}
  \left| \int (e^{itg}-1) w''_t \right| \leq 2 \int_{\omega \geq |t|^{-1/10}}
  |w''_t|
  \leq 2 \int_{\omega \geq |t|^{-1/10}} (1+C |t|^{1/10})^ {\omega}
  \leq C\sum_{n=|t|^{-1/10}}^\infty \rho^n (1+C |t|^{1/10})^n.
  \end{equation*}
When $t$ is small enough, $\rho (1+C |t|^{1/10}) <\sqrt{\rho}< 1$.
Hence, $  | \int (e^{itg}-1) w''_t | \leq C \rho^{|t|^{-1/10}/2} =
o(t^2 L(1/|t|))$. In the same way, $\left | \int g w''_t \right|
\leq \norm{g}_{L^{3/2}} \norm{w''_t}_{L^3}$ and
$\norm{w''_t}_{L^3}$ decays stretched exponentially, whence it is
$o(t^2 L(1/|t|))$. Finally,
  \begin{equation*}
  \left| \int 1_{|g| \leq 1/|t|}
  g^2 w''_t \right| \leq \frac{1}{|t|^2} \int |w''_t| =O(\rho^{
|t|^{-1/10}/2} / t^2) = o(t^2 L(1/|t|)).
  \end{equation*}
Hence, \eqref{sur_w't} holds also for $w''_t$. We get
  \begin{equation*}
  \int (e^{itg}-1-itg) w_t=  -\frac{t^2}{2} \int 1_{|g| \leq 1/|t|}
  g^2 w_t+o(t^2 L(1/|t|)).
  \end{equation*}
Since $\int e^{itg} w_t =\lambda_t$ and $\int w_t=1$, this gives
the conclusion of the lemma.
\end{proof}

\begin{cor}
\label{estime_lambdat}
We have $\lambda_t=1+O(|t|^{11/10})$.
\end{cor}
\begin{proof}
In the proof of the previous lemma, we have proved that $\int_{
\omega>|t|^{-1/10}} 1_{|g|\leq 1/|t|} g^2 w_t=O(t)$. Moreover, on
$\{ \omega\leq  |t|^{-1/10}\}$, the function $w_t$ is uniformly
bounded. Hence, $\int_{ \omega\leq |t|^{-1/10}} 1_{|g|\leq 1/|t|}
g^2 w_t \leq C \int 1_{|g|\leq 1/|t|} g^2 \sim C L(1/|t|)$. Hence,
  \begin{equation*}
  \int 1_{|g| \leq 1/|t|} g^2 w_t = O(L(1/|t|)).
  \end{equation*}
Moreover, by \eqref{w_t_proche_1} and our assumption $\int
g=0$,
  \begin{equation*}
  \left| \int gw_t \right|
  =\left| \int g(w_t-1) \right|
  \leq \norm{g}_{L^{3/2}} \norm{w_t-1}_{L^3} =
  O(|t|^{1/10}).
  \end{equation*}
This proves that, in \eqref{decrit_lambdat_1}, the right side is
$1+O(|t|^{11/10})$.
\end{proof}

\subsubsection{Estimates on the basis}

To proceed, we will need to know that $w_t$ is constant on the
basis up to $O(t)$. We already know that
$\norm{w_t-1}_{L^3}=O(|t|^{1/10})$, but this is not sufficient to
estimate precisely the terms in \eqref{decrit_lambdat_1}. To get
such an estimate, we will need real continuity, and not only the
weak continuity given by Keller-Liverani's theorem. This will be
achieved by working directly on the basis. The goal of this
paragraph is to prove Lemma~\ref{almost_ct_on_0}.

Let $U_0$ be the map induced by $U$ on the basis $\Delta_0$ of the
tower. Denote by $\phi$ the first return time to the basis, so
that $U_0(x)=U^{\phi(x)}(x)$.

Let us consider the  space $\boH_0$ of H\"older functions
$u:\Delta_0\to \C$ on the basis, and define an operator
$R_n:\boH_0\to \boH_0$ by $R_n u(x) =\sum J^{(n)}(y) u(y)$, where
the sum is restricted to those $y\in \Delta_0$ with return time
$\phi(y)=n$, and $U^n(y)=x$. Set also $R_n(t) (u)=R_n(e^{it S_n g}
u)$.

\begin{lem}
There exist $C>0$ and $\theta<1$ such that, for all $n \in \N$ and
all $t\in [-1,1]$, $\norm{R_n(t)} \leq C \theta^n$ and
$\norm{R_n(t)-R_n} \leq C \theta^n |t|$.
\end{lem}
\begin{proof}
The map $U_0$ is Gibbs-Markov on $\Delta_0$. Hence, \cite[Lemma
3.2]{gouezel:stable} proves that $\norm{R_n} \leq C \mu_\Delta(
\phi=n)$ and \cite[Lemma 3.5]{gouezel:stable} yields
$\norm{R_n(t)-R_n} \leq C|t| n \mu_\Delta(\phi=n) + C
\int_{\phi=n} |e^{itS_n g}-1|$.

Since $\mu_\Delta(\phi=n)=O(\rho^n)$, we get in particular
$\norm{R_n} \leq C \rho^n$, which decays exponentially. Moreover, on
$\{\phi=n\}$, $|S_n g|^{3/2} \leq n^{1/2} \sum_{k=0}^{n-1} |g\circ U^k
|^{3/2}$, whence $\int_{\phi=n} |S_n g|^{3/2} \leq n^{1/2}
\int_{\Delta} |g|^{3/2} = O(n^{1/2})$. Consequently,
  \begin{equation*}
  \norm{R_n(t)-R_n} \leq C|t| n \rho^n + C \int 1_{\phi=n} |t| |S_n g|
  \leq C|t| n\rho^n +C |t| \norm{1_{\phi=n}S_n g}_{L^{3/2}} \norm{
  1_{\phi=n}}_{L^3},
  \end{equation*}
which decays also exponentially.
\end{proof}

For $|z|<\theta^{-1}$, it is possible to define $R(z,t):= \sum z^n
R_n(t)$. The operator $R(1,0)$ is the transfer operator associated
to $U_0$. It has a simple isolated eigenvalue at $1$, and the
corresponding eigenfunction is the constant function $1$. Hence,
for $(z,t)$ close enough to $(1,0)$, $R(z,t)$ has a unique
eigenvalue $\lambda(z,t)$ close to $1$.

\begin{lem}
\label{lem:Rzt}
We have $\norm{ R(z,t) -R(1,0)} =O( |t|+|z-1|)$.
\end{lem}
\begin{proof}
We have $R(z,t)-R(1,t)= \sum (z^n -1) R_n(t)$. Moreover, if
$|z|\leq \theta^{-1/2}$, $|z^n -1|\leq |z-1| \sum_{k=0}^{n-1}
|z|^k \leq C |z-1| \theta^{-n/2}$. Hence,
  \begin{equation*}
  \norm{ R(z,t)-R(1,t)} \leq \sum_{n=0}^\infty |z^n-1| \norm{R_n(t)}
  \leq C |z-1| \sum_{n=0}^\infty  \theta^{-n/2} \theta^n \leq
  \frac{C}{1-\theta^{1/2}} |z-1|.
  \end{equation*}
Moreover,
  \begin{equation*}
  \norm{R(1,t)-R(1,0) } \leq \sum_{n=0}^\infty \norm{R_n(t)-R_n}
  \leq \sum_{n=0}^\infty C|t| \theta^n \leq \frac{C}{1-\theta} |t|.
  \qedhere
  \end{equation*}
\end{proof}

\begin{lem}
\label{eigenvalue_restricted}
We have $R(\lambda_t^{-1}, t) (1_{\Delta_0} w_t)=1_{\Delta_0}
w_t$.
\end{lem}
\begin{proof}
Let $x\in \Delta_0$, let $\{x_p\}$ be the set of its preimages
under $U$, at respective heights $\omega_p$, and let $x'_p$ be the
projection of $x_p$ in the basis. Since $\hat{U}_t w_t=\lambda_t
w_t$, we have $\lambda_t w_t(x)=\sum e^{it g(x_p)} J(x_p)
w_t(x_p)$. Moreover, for any $y\in \Delta$ with $\haut(y)>0$, we
have $\lambda_t w_t(y)=e^{itg(U^{-1}y)} w_t(U^{-1}y)$. Hence,
$\lambda_t^{\omega_p}w_t(x_p)=e^{it S_{\omega_p} g(x'_p)}
w_t(x'_p)$. We get
  \begin{equation}
  \label{in_the_basis}
  w_t(x) =\sum \lambda_t^{-\omega_p-1} J^{(\omega_p+1)}(x'_p) e^{it S_{\omega_p+1}
  g(x'_p)}  w_t(x'_p).
  \end{equation}
The points $x'_p$ are exactly the preimages of $x$ under $U_0$,
and the corresponding return time for $U$ is $\omega_p+1$. Hence,
\eqref{in_the_basis} gives the conclusion of the lemma.
\end{proof}

We have all the necessary tools to prove the main result of this
paragraph:
\begin{lem}
\label{almost_ct_on_0}
For $t\in [-1,1]$, there exists $c(t)\in \C$ such that
$\norm{1_{\Delta_0} (w_t-c(t))}_\infty = O(t)$. Moreover, $c(t)\to
1$ when $t\to 0$.
\end{lem}

\begin{proof}
Lemma~\ref{eigenvalue_restricted} proves that
$\lambda(\lambda_t^{-1}, t)=1$, and the corresponding
eigenfunction is $1_{\Delta_0}w_t$. Let $Q_t$ be the
eigenprojection of $R(\lambda_t^{-1},t)$ corresponding to the
eigenvalue $1$. It satisfies
  \begin{equation*}
  \norm{Q_t -Q_0} = O(\norm{R(\lambda_t^{-1},t)-R(1,0)}) =
  O( |\lambda_t^{-1}-1| +|t|)
  =O(|t|)
  \end{equation*}
by Lemma~\ref{lem:Rzt} and Corollary~\ref{estime_lambdat}. Let
$b_t = Q_t 1_{\Delta_0}$. As $b_0=1_{\Delta_0}$, $b_t$ satisfies
$\norm{b_t-1_{\Delta_0}}=O(t)$. In particular, $b_t \to
1_{\Delta_0}$ in $L^1$.

The function $b_t$ is proportional to $w_t$ on the basis
$\Delta_0$. Hence, there exists a scalar $c(t)$ such that
$1_{\Delta_0}w_t=c(t) b_t$. Since $w_t$ goes to $1$ in $L^1$ when
$t\to 0$, we get
  \begin{equation*}
  c(t)=\frac{\int 1_{\Delta_0} w_t}{\int b_t} \to \frac{\mu_\Delta(\Delta_0)}{\int
  b_0}=1.
  \end{equation*}
Finally,
  \begin{equation*}
  \norm{1_{\Delta_0}(w_t-c(t))}_\infty=|c(t)| \norm{ b_t
  -b_0}_\infty = O(t).
  \qedhere
  \end{equation*}
\end{proof}

\subsubsection{Proof of Theorem~\ref{thm_ad_young} when $\int g=0$}

\label{cas_particulier}

Let the function $G$ be given by $G(x)=\sum_{k=0}^{\haut(x)-1}
g(U^k \pi_0 x)$, as in Theorem~\ref{thm_ad_young}.

\begin{lem}
\label{lemme_estime_lambdat_2}
When $t\to 0$,
  \begin{equation}
  \label{decrit_lambdat_2}
  \lambda_t = 1-(1+o(1))\frac{t^2}{2}\int 1_{|g|\leq 1/|t|} g^2 e^{itG}
 + it(1+o(1)) \int g e^{itG} + o(t^2 L(1/|t|)).
  \end{equation}
\end{lem}
\begin{proof}
We will start from \eqref{decrit_lambdat_1} and show that we can
replace $w_t$ by $e^{itG}$.

We have $w_t(x)=\lambda_t^{-\haut(x)} e^{it G(x)} w_t(\pi_0 x)$.
Hence, by Corollary~\ref{estime_lambdat} and
Lemma~\ref{almost_ct_on_0}
  \begin{align*}
  | w_t(x)-c(t) e^{itG(x)}| &= |\lambda_t^{-\haut(x)}w_t (\pi_0 x) -c(t)|
  \leq |\lambda_t^{-\haut(x)}-1| | w_t(\pi_0 x)| + |w_t(\pi_0 x)-c(t)|
  \\&
  \leq \bigl[(1+C|t|^{11/10})^{\haut(x)}-1 \bigr] C +C|t|
  \leq \haut(x)C|t|^{11/10} (1+C|t|^{11/10})^{\haut(x)} +C|t|.
  \end{align*}
Fix $b>0$ large enough. For $\haut(x) \leq b \log (1/|t|)$, we
obtain $|w_t(x)- c(t)e^{itG(x)}| \leq C|t|$. For $\haut(x) \geq b
\log(1/|t|)$ and small enough $t$, we also get
$|w_t(x)-c(t)e^{itG(x)}| \leq \rho^{-\haut(x)/4}$.

Hence,
  \begin{equation*}
  \int_{\omega\geq b \log(1/|t|)} 1_{|g|\leq 1/|t|}g^2 |w_t -c(t)e^{itG}|
  \leq \int_{\omega\geq b \log(1/|t|)} \frac{1}{|t|^2} \rho^{-\haut(x)/4}
  \leq \frac{1}{|t|^2}C \sum_{n=b\log (1/|t|)}^\infty \rho^n \rho^{-n/4}
  =o(1)
  \end{equation*}
if $b$ is large enough. Moreover,
  \begin{equation*}
  \int_{\omega\leq b \log(1/|t|)} 1_{|g|\leq 1/|t|}g^2 |w_t -c(t)e^{itG}|
  \leq C|t| \int 1_{|g| \leq 1/|t|} g^2 = C|t| L(1/|t|).
  \end{equation*}
Hence,
  \begin{equation}
  \label{kjahsdljhaklsdf}
  \int 1_{|g| \leq 1/|t|} g^2 w_t =c(t) \int 1_{|g| \leq 1/|t|} g^2
  e^{itG} +o(1)
  = (1+o(1)) \int 1_{|g| \leq 1/|t|} g^2e^{itG} + o(1).
  \end{equation}

In the same way,
  \begin{equation*}
  \int_{\omega \geq b \log(1/|t|)} |g| |w_t-c(t)e^{itG}|
  \leq \norm{g}_{3/2} \norm{ 1_{\omega\geq b \log(1/|t|)}
  |w_t-c(t)e^{itG}|}_{L^3}
  =O(t)
  \end{equation*}
if $b$ is large enough. Moreover,
  \begin{equation*}
  \int_{\omega \leq b \log(1/|t|)} |g| |w_t-c(t)e^{itG}|
  \leq \int |g| C|t|=O(t).
  \end{equation*}
We get
  \begin{equation}
  \label{lkafkljhksladf}
  \int g w_t = c(t) \int g e^{itG} +O(t)=(1+o(1))\int g e^{itG} +O(t).
  \end{equation}
Equations \eqref{kjahsdljhaklsdf} and \eqref{lkafkljhksladf}
together with \eqref{decrit_lambdat_1} imply
\eqref{decrit_lambdat_2}.
\end{proof}

\begin{rmq}
The proof of the lemma also shows that, in
\eqref{decrit_lambdat_2}, it is sufficient to integrate on $\{
\haut(x) \leq b \log(1/|t|)\}$ if $b$ is large enough, since the
remaining part is in $o(t^2 L(1/|t|))$.
\end{rmq}

The following lemma will use the additional assumptions that $l$
is slowly varying and $l(x\ln x) \sim l(x)$, $L(x\ln x) \sim
L(x)$.
\begin{lem}
We have
  \begin{equation}
  \label{woyrtuiert}
  \int 1_{|g|\leq 1/|t|} g^2 e^{itG} = L(1/|t|) (1+o(1)).
  \end{equation}
\end{lem}
\begin{proof}
It is sufficient to prove \eqref{woyrtuiert} on $\{ \omega\leq b
\log(1/|t|)\}$, since the remaining part can be ignored.

Take some $\epsilon>0$, we will prove that
  \begin{equation}
  \label{qwoiurtwqwr}
  \int_{\omega \leq b\log(1/|t|)} 1_{|g|\leq 1/|t|} g^2 |e^{itG}-1|
  \leq 2\epsilon L(1/|t|)
  \end{equation}
when $t$ is small enough. This will be sufficient to conclude the
proof.

Let $A_t:=\{ x\tq \haut(x) \leq b\log (1/|t|), |G(x)|\geq
\epsilon/|t|\}$. if $x\in A_t$, there exists $y$ below $x$ in the
tower such that $|g(y)|\geq \frac{\epsilon}{|t| b\log(1/|t|)}$.
Let $B=\{ x \tq |g(x)|\geq \frac{\epsilon}{|t| b \log(1/|t|)}\}$,
we get $\mu_\Delta(A_t) \leq b \log(1/|t|) \mu_\Delta(B)$.

Let $Z$ be a random variable on $\R$ with the distribution of $g$.
Then
  \begin{align*}
  P\left( \frac{1}{|t| \log(1/|t|)^2} \leq |Z| \leq 1/|t|\right)&
  = |t|^2 \log(1/|t|)^4 l\left( \frac{1}{|t| \log(1/|t|)^2}\right)
  -|t|^2 l(1/|t|)
  \\&
  = |t|^2 \log(1/|t|)^4 l(1/|t|) (1+o(1))
  \end{align*}
and
  \begin{align*}
  b \log(1/|t|) P\left( |Z| \geq \frac{\epsilon}{|t| b
  \log(1/|t|)}\right)
  &
  = b \log(1/|t|) \frac{ |t|^2 b^2 \log(1/|t|)^2}{\epsilon^2} l\left(
  \frac{\epsilon}{|t| b \log(1/|t|)}\right)
  \\&
  = \frac{|t|^2 b^3 \log(1/|t|)^3}{\epsilon^2}l(1/|t|) (1+o(1)).
  \end{align*}

Hence, if $t$ is small enough, $\mu_\Delta(A_t) \leq P\left(
\frac{1}{|t| \log(1/|t|)^2} \leq |Z| \leq 1/|t|\right)$.

We would like to estimate $\int_{A_t} 1_{|g| \leq 1/|t|} g^2$. Now
\begin{equation*}
\int_{A_t} 1_{|g| \leq 1/|t|} g^2= \int_{A_t} 1_{ \frac{1}{|t|
\log(1/|t|)^2}< |g| \leq 1/|t|} g^2 + \int_{A_t} 1_{ |g| \leq
\frac{1}{|t| \log(1/|t|)^2}} g^2.
\end{equation*}
On the one hand
\begin{equation*}
\int_{A_t} 1_{ \frac{1}{|t| \log(1/|t|)^2}< |g| \leq 1/|t|} g^2
\leq \int_{\frac{1}{|t| \log(1/|t|)^2}}^{1/|t|} x^2 \dd P(x),
\end{equation*}
and on the other hand, by applying the above bounds we get
\begin{align*}
\int_{A_t} 1_{ |g| \leq \frac{1}{|t| \log(1/|t|)^2}} g^2 &\leq
\frac{1}{|t|^2 \log(1/|t|)^4} \mu_{\Delta}(A_t) \\& \leq
\frac{1}{|t|^2 \log(1/|t|)^4} P\left( \frac{1}{|t| \log(1/|t|)^2}
\leq |Z| \leq 1/|t|\right) \\& \leq \int_{\frac{1}{|t|
\log(1/|t|)^2}}^{1/|t|} x^2 \dd P(x).
\end{align*}

Thus we need to deal with the integral
  \begin{equation*}
  \int_{\frac{1}{|t|
\log(1/|t|)^2}}^{1/|t|} x^2 \dd P(x),
  \end{equation*}
which is equal to
  \begin{equation*}
  L(1/|t|)-L\left( \frac{1}{|t|
  \log(1/|t|)^2} \right) = L(1/|t|) - L(1/|t|)(1+o(1))=o(L(1/|t|)).
  \end{equation*}

Hence, for small enough $t$, we get
  \begin{equation*}
  \int_{A_t} 1_{|g| \leq 1/|t|} g^2 | e^{itG}-1| \leq \epsilon L(1/|t|).
  \end{equation*}

On $B_t:=\{ x\tq \haut(x)\leq b\log(1/|t|), x\not \in A_t\}$, we
have $|e^{it G(x)}-1| \leq |t| |G(x)|\leq \epsilon$. Hence,
  \begin{equation*}
  \int_{B_t} 1_{|g|\leq 1/|t|} g^2 |e^{itG}-1| \leq \epsilon \int
  1_{|g|\leq 1/|t|} g^2 =\epsilon L(1/|t|).
  \end{equation*}
These two equations imply \eqref{qwoiurtwqwr}. This concludes the
proof.
\end{proof}

Since $\int g=0$, Lemma~\ref{lemme_estime_lambdat_2} gives
  \begin{align*}
  \lambda_t&=1-\frac{t^2}{2}L(1/|t|)(1+o(1))+ (1+o(1))it\int g e^{itG}
  \\&
  = 1-\frac{t^2}{2}L(1/|t|)(1+o(1))+ (1+o(1))it\int g (e^{itG}-1)
  \\&
  =1- \frac{t^2}{2}L_1(1/|t|) (1+o(1)),
  \end{align*}
since $\int g (e^{itG}-1) = i(a+o(1))t L(1/|t|)$ by assumption.

This asymptotic expansion readily implies the conclusion of
Theorem~\ref{thm_ad_young}, for $g$ such that $\int g=0$.

\subsubsection{Proof of Theorem~\ref{thm_ad_young} in the general case}

\label{etend_cas_general}

Let $g:\Delta \to \R$ belong to $L^p(\Delta)$ for any $p<2$ (this
is in particular the case if $g$ satisfies the assumptions of
Theorem~\ref{thm_ad_young}). Set $G(x)=\sum_{k=0}^{\haut(x)-1}
g(U^k \pi_0 x)$.

\begin{lem}
\label{integrable}
For any $p<2$, the function $G$ belongs to $L^p(\Delta)$.
\end{lem}
\begin{proof}
For $x\in \Delta$, let $\phi(x)$ be its return time to the basis.
Set also $\psi(x)= \phi(\pi_0 x)$, where $\pi_0$ is the projection
on the basis of the tower.

We have
  \begin{equation*}
  |G(x)|^p = \left| \sum_{k=0}^{\haut(x)-1} g(U^k \pi_0 x) \right|^p
  \leq \haut(x)^{p-1} \sum_{k=0}^{\haut(x)-1} |g(U^k \pi_0 x)|^p.
  \end{equation*}
Changing variables, we get
  \begin{equation*}
  \int |G(x)|^p \leq \int |g(y)|^p \sum_{k=1}^{\phi(y)-1} \haut(U^k
  y)^{p-1}
  \leq \int |g(y)|^p \psi(y)^p.
  \end{equation*}
Since the tower has exponentially small tails, the function $\psi$
belongs to $L^q$ for any $q<\infty$. Using the H\"older inequality
with a sufficiently large $q$, we obtain $\int |G(x)|^p <\infty$.
\end{proof}

Let $g'$ be another function on $\Delta$. Define also $G'(x)=
\sum_{k=0}^{\haut(x)-1} g'(U^k \pi_0 x)$.
\begin{lem}
\label{lem_equiv_bar}
if $g-g'$ is bounded, then
  \begin{equation*}
  \int g (e^{itG}-1) = \int g'(e^{it G'}-1) + O(t)
  \end{equation*}
when $t\to 0$.
\end{lem}
\begin{proof}
We have
  \begin{equation*}
  \int g(e^{itG}-1) - \int  g' (e^{it  G'}-1)
  =\int (g- g') (e^{itG'}-1) + \int  g (e^{itG}- e^{it  G'}).
  \end{equation*}
Since $g- g'$ is bounded, the first integral satisfies
  \begin{equation*}
  \left| \int (g- g') (e^{itG'}-1) \right| \leq C|t| \int |G'|,
  \end{equation*}
which is $O(t)$ since $G'$ is integrable by Lemma
\ref{integrable}. For the second integral, $|G(x)- G'(x)|\leq C
\haut(x)$. Hence,
  \begin{equation*}
  \left|\int  g (e^{itG}- e^{it  G'}) \right|
  \leq \int |g| C|t| \omega \leq C|t| \norm{g}_{L^{3/2}}
  \norm{\omega}_{L^3}
  =O(t).
  \qedhere
  \end{equation*}
\end{proof}

\begin{proof}[Proof of Theorem~\ref{thm_ad_young}]
Let $g$ satisfy the assumptions of Theorem~\ref{thm_ad_young}.
Write $ g'= g -\int g$. Then $ g'$ is still in the nonstandard
domain of attraction of the normal law, and its distribution
functions $ l'$ and $ L'$ satisfy $ l' \sim l$ and $ L' \sim L$.
Moreover,
  \begin{equation*}
  \int  g' (e^{it G'}-1) = (a+o(1))it L(1/|t|)),
  \end{equation*}
since $g$ satisfies the same estimate and Lemma
\ref{lem_equiv_bar} applies.

We have already proved Theorem~\ref{thm_ad_young} for functions of
zero integral. This applies to $ g'$, and gives $\frac{S_n  g'
}{B_n} \to \boN(0,1)$. Since $S_n  g'= S_n g-n\int g$, this
concludes the proof.
\end{proof}

\begin{rmq}
Lemmas~\ref{integrable} and \ref{lem_equiv_bar} do not involve the
dynamics of the returns to the basis. Hence, they also hold in
hyperbolic Young towers.
\end{rmq}

\section{Estimate of the integral in the stadium billiard}

Let us turn back to the study of the stadium. We will use the
notations of the first two sections. In particular, starting from a
fixed function $f_0:X_0 \to \R$ satisfying $(P1)$, we have
obtained a function $g:\Delta \to \R$.
According to Theorem~\ref{thm_ad_young}, if we want to obtain a
limit theorem for $g$, the quantity to be estimated is $\int g
(e^{itG}-1)$. The main result of this section is the following
proposition.

\begin{prop}
\label{estime_integrale_precise} Let $y=\frac{1}{1-\frac{3}{4}\log
3}$, and recall the definition of $I$ from \eqref{donne_I}. We
have
  \begin{equation*}
  \int_{\Delta} g (e^{itG}-1)\dd\mu_\Delta=i\frac{I^2
  (y-1)\ell^2}{\pi}t \log(1/|t|) + o(t\log(1/|t|)).
  \end{equation*}
\end{prop}

Our proof approximates the left hand side with an integral
explicitly given in the phase space of the stadium. This later
integral can be estimated with sufficient precision due to the
geometric properties of the billiard map.

\subsection{Preliminary estimates}

First we show that the relevant expression can be pulled back to
the hyperbolic Young tower. Let $\omega(x)$ be the height of the
point $x$ in $\bar\Delta$, and let $\bar\pi_0:\bar\Delta \to
\bar\Delta_0$ be the projection on the basis. We define two
functions $\bar F$ and $\bar G$ on $\bar \Delta$ by $\bar
F(x)=\sum_{k=0}^{\omega(x)-1} \bar f(\bar U^k \bar \pi_0 x)$ and
$\bar G(x)=\sum_{k=0}^{\omega(x)-1} \bar g(\bar U^k \bar \pi_0
x)$.

\begin{lem}
\label{lem_tower1}
We have
  \begin{equation*}
  \int_{\Delta} g (e^{itG}-1)=\int_{\bar \Delta} \bar f (e^{it \bar
  F}-1)+ O(t).
 \end{equation*}
\end{lem}
\begin{proof}
As $(\pi_\Delta)_*(\mu_{\bar \Delta})=\mu_\Delta$, $\int_\Delta g
(e^{it G}-1)=\int_{\bar \Delta} \bar g (e^{it \bar G}-1)$
automatically. As $\bar g-\bar f$ is bounded by
Lemma~\ref{u_bornee}, Lemma~\ref{lem_equiv_bar} implies the
statement.
\end{proof}

Note that $\bar F$ is essentially a Birkhoff sum of $f$ for the
\emph{inverse map} $T^{-1}$. Thus, if we switch from $T$ to
$T^{-1}$, we may investigate our integral by dynamical tools.

For all $x\in X$, let $\gam(x)=f(T^{-1}x)$. Introduce $\bar
\gam=\gam \circ \pi_X$. For $x\in \bar\Delta$ with $\haut (x)>0$,
let $\bar \Gam(x)= \sum_{k=1}^{\haut(x)-1} \bar \gam (\bar
U^{-k}x)$, or equivalently, $\bar \Gam(x) =\bar F(\bar U^{-1}x)$.
We fix $\bar \Gam(x)=0$ on $\Delta_0$.

\begin{lem}
\label{lem_tower2}
We have
  \begin{equation*}
  \int \bar f (e^{it\bar F}-1) =
  \int \bar \gam
  (e^{it\bar \Gam}-1) +O(t).
  \end{equation*}
\end{lem}
\begin{proof}
 We have $\bar \gam \circ \bar
U=\bar f$, and, apart from $\bar U^{-1}(\bar \Delta_0)$, $\bar
\Gam \circ \bar U=\bar F$. Thus,
  \begin{align*}
  \int \bar \gam (e^{it\bar \Gam}-1)-
  \int \bar f (e^{it\bar F}-1)
  &
  =
  \int \bar \gam\circ \bar U
  (e^{it\bar \Gam\circ \bar U}-1)-
  \int \bar f (e^{it\bar F}-1)
  \\&
  =  \int_{\bar U^{-1}(\bar \Delta_0)}\bigl[\bar \gam\circ \bar U
  (e^{it\bar \Gam\circ \bar U}-1)-
  \bar f (e^{it\bar F}-1)\bigr].
  \end{align*}
As $\phi_+$ is bounded on the rectangle $R$ that defines the basis
of the tower, the functions $\bar \gam \circ \bar U$ and $\bar f$
are bounded on $\bar U^{-1}(\bar \Delta_0)$. By
Lemma~\ref{integrable} $\bar F$ and $\bar \Gam$ are both
integrable. This completes the proof.
\end{proof}

We will consider $T^{-1}$ as the first return map of $T_0^{-1}$ to
the subspace $X$. The return time is $\phi_-=\phi_+ \circ T^{-1}$.

\subsection{Geometric properties of $T^{-1}$ in the vicinity of
its singularities}

The behavior of $\int \bar \gam(e^{it\bar \Gam}-1)$ is governed by
the dynamical properties of $T^{-1}$ at those parts of the phase
space where it is equivalent to a long series of bounces along the
parallel segments of the boundary. These sets have the following
structure: the points for which $T^{-1}$ acts as $n$ consecutive
bounces on the segments form two stripes of slope approximately
$-1$. $T^{-1}$ maps these two cells onto two stripes of positive
slope. The geometry is described on the figure below.

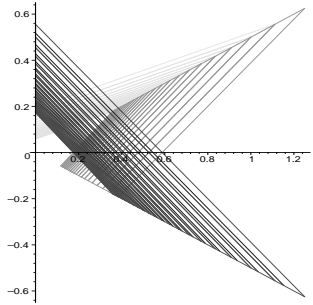
\begin{center}
\includegraphics[width=4cm,height=4cm]{sing.eps}
\end{center}

More precisely, this figure appears four times in $X$ (twice in
both of the parallelograms that define $X$). The transformation
$T^{-1}$ jumps from one such region to another, however, they play
the same dynamical role. Thus, to simplify matters, we pretend as
if we had only one of them.

Recall that $\phi_-=\phi_+\circ T^{-1}$ is the return time of
$T_0^{-1}$ to $X$. As a further notational simplification we
ignore that there are two stripes on which $\{\phi_-=n\}$. Let
$M_n$ stand for the stripe $\{\phi_-=n\}$, which will be also
referred to as the set of points of \emph{return time} $n$.

Later on we shall see that the other type of singularities
(corresponding to trajectories sliding along the semi-circle) does
not play any role in the leading term behavior of $\int \bar \gam
(e^{it\bar \Gam}-1)$.

\begin{rmq}
We need to study the map $T^{-1}$ and not $T$. These two are not
isomorphic, as $X$ is the set of points on a semi-circle for which
the \emph{previous} collision is not on that semi-circle. This
definition introduces an asymmetry of past and future. The map
$T^{-1}$ is, however,  isomorphic to the map  induced on the set
of points on a semi-circle for which the \emph{next} collision is
on another semi-circle. This later induced map has been studied by
Markarian in \cite{markarian:slow}, where he, in particular, has
shown that it satisfies Chernov's axioms from
\cite{chernov:decay}.
\end{rmq}

Fix $\rho<1$ such that the tails of the original Young tower are
bounded by  $c \rho^n$, and $K>0$ such that $K\log(\rho^{-1})>4$.
In what follows we essentially consider curves with tangent
vectors in the unstable cone of $T^{-1}$ (u-curves of $T^{-1}$ for
short). 
\begin{comment}
Let $C_D$ be a constant that bounds the distortions of all
iterates of $T^{-1}$ along u-curves (as long as the iterated
points belong to the same cells).
\end{comment}

\begin{defn}
Consider the stripe $M_n$ and its two sides of slope $-1$. A
\emph{good curve} $\boC$ of return time $n$ is a $C^1$ curve that
connects these two sides. We put further requirements on the slope
of $\boC$: it should belong to the interval $[1/4,4]$ for all
points and it should be constant up to $1/\sqrt{n}$ precision
(i.e., for all points $x$ and $y$ in $\boC$, the slopes of $\boC$
at $x$ and $y$, $s(x)$ and $s(y)$, should satisfy $|s(x)-s(y)|\leq
\frac{1}{\sqrt{n}}$).
\end{defn}

\begin{rmq}
\label{good_even_better} Note that our requirement on the slope in
\emph{not} a curvature bound. Stated in this form, it is not hard
to see that good curves tend to have more and more constant slopes
when iterated by $T^{-1}$. To see this consider a good curve with
large return time $n$ and iterate it backwards by the billiard
flow. Just before collision it corresponds to a dispersing
wavefront that defocuses within finite time and, while
experiencing many bounces with the straight walls, loses most of
its curvature. At the time moment just after the previous
collision on the other semi-circle, this wavefront is flat up to
$1/n$. Thus, any subcurve $\boC'\subset T^{-1}\boC$ is
automatically a good curve if it connects the two long sides of a
stripe.
\end{rmq}

\begin{defn}
A \emph{standard} curve is defined as a good curve of constant
slope $1$. In particular, it is a segment.
\end{defn}
The choice of $1$ as the slope for standard curves is arbitrary.
More important is the fact that the standard curves of return time
$n$ give a fixed foliation for (most of) the stripe $M_n$.

If $\boC$ is a good curve of return time $n$, any point of
$T^{-1}(\boC)$ has return time at least $n/3-\Cmm$ and at most
$3n+\Cmm$ for some constant $\Cmm$. Furthermore, there exists a
constant $\Cn\in \N$ such that for any $i\in [n/3+\Cn, 3n-\Cn]$,
the set $T^{-1}(\boC) \cap \{\phi_-=i\}$ is a good curve of return
time $i$ (see Remark~\ref{good_even_better}). Let us denote
$\boC_i=\{x\in \boC \tq \phi_-(T^{-1}x)=i\}$. We also have
  \begin{equation}
  \label{mesure_jete}
  \frac{\Leb \left( \boC \moins \bigcup_{i\in [n/3+\Cn, 3n-\Cn]}
  \boC_i\right)} {\Leb(\boC)} \leq \frac{\Ca}{n}
  \end{equation}
for a universal constant $\Ca$. We will say that the set $\boC
\moins \bigcup_{i\in [n/3+\Cn, 3n-\Cn]} \boC_i$ is \emph{thrown
away} at the first iterate of $\boC$. Formula \eqref{mesure_jete}
shows that the points which are thrown away have negligible
measure.

\begin{rmq}
\label{transition_proba} In addition to the above observations, it
is possible to estimate the transition probabilities from one
stripe to the other in the following sense. There exists a
sequence $\epsilon_n$ that tends to $0$ as $n \to \infty$, such
that for any good curve ${\boC}$ of return time $n$, and for any
$i\in [n/3+\Cn,3n-\Cn]$,
  \begin{equation}
  \label{proba_transition}
  (1-\epsilon_n) \frac{3n}{8i^2}\leq
  \frac{\Leb\{ x\in {\boC} \tq \phi_-(T^{-1}(x))=i\}}{\Leb({\boC})}
  \leq \frac{3n}{8i^2}
  (1+\epsilon_n).
  \end{equation}
This can be verified by direct calculation. In other words, we go
from $n$ to $i$ asymptotically with probability $\frac{3n}{8i^2}$
(note that $\sum_{i=n/3+\Cn}^{3n-\Cn} \frac{3n}{8i^2} \to 1$).
\end{rmq}

Applying the above process several times, we may iterate the good
curves by $T^{-1}$ and obtain finer and finer partitions of
$\boC$. A sequence of integers $n_0,n_1,\dots,n_k$ is referred to
as \emph{admissible} if, for all $i<k$, $n_{i+1} \in [n_i/3+\Cn,
3n_i-\Cn]$. Given a good curve of return time $n_0$, $\boC$, and
an admissible sequence $n_0,\dots,n_k$, let
  \begin{equation*}
  \boC_{n_0,\dots,n_k}=\{ x\in \boC \tq \forall i \leq k,
  \phi_-(T^{-i}x)=n_i\}.
  \end{equation*}
This is a subcurve of $\boC$ mapped by $T^{-k}$ onto a good curve
of return time $n_k$.

\begin{lem}
\label{comparables} There exists a constant $\Co>0$ such that, for
any pair of good curves of the same return time $n_0$, $\boC$ and
$\boC'$, and for any fixed admissible sequence $n_0,\dots,n_k$, we
have
  \begin{equation*}
  \Co^{-1} \leq \frac{
  \Leb(\boC_{n_0,\dots,n_k})}{\Leb(\boC'_{n_0,\dots,n_k})} \leq \Co.
  \end{equation*}
\end{lem}
\begin{proof}
This follows from the uniform expansion and the bounded distortion
properties of $T^{-1}$ along its u-curves.
\end{proof}

In what follows, when we talk about iterating a good curve, we
will always mean the above process of refinement, along with
throwing away some part at each step. However, the number of
iterations may depend on the point of $\boC$ we are considering.
This is formulated in the following definition.

\begin{defn}
Let $\boC$ be a good curve of return time $n$. Let furthermore $A$
be a subset of $\boC$ and $\tau:\boC\moins A \to \N$. Then
$(A,\tau)$ is a \emph{stopping time} on $\boC$ if
\begin{itemize}
\item There exists $p\in \N$ such that $3^{p+1}<n_0$, 
with the following property: all the
connected components of $\boC\moins A$ are of the form
$\boC_{n_0,\dots,n_k}$, where $n_0=n$, the sequence
$n_0,\dots,n_k$ is admissible, and $n_k \in [3^p, 3^{p+1} -1]$.
Furthermore, $\tau$ is uniformly equal to $k$ on such a component.
\item We have $\Leb(A)/\Leb(\boC)\leq 1/2$.
\end{itemize}
\end{defn}
Here typically $3^p\ll n$, thus we stop at the first occasion when
the return time decreases below a certain level.

\begin{rmq}
If $(A,\tau)$ is a stopping time on $\boC$, then
  \begin{equation*}
  \frac{1}{2} \leq
  \frac{\Leb(\boC\moins A)}{\Leb(\boC)} \leq 1.
  \end{equation*}
Thus in our estimates $\Leb(\boC \moins A)$ and  $\Leb(\boC)$ may
be replaced with each other. We will often use this without giving
further details.
\end{rmq}

Let us define, in particular, the \emph{standard stopping time}
for a good curve $\boC$ of return time $n$. Let $p$ be the integer
for which $3^p \leq n^{1/4} < 3^{p+1}$. We partition $\boC$,
iterate $T^{-1}$ and throw away the negligible parts according to
the process described above. We go on iterating until either the
return time of the image belongs to the interval
$[3^p,3^{p+1}-1]$, or the number of iterates exceeds $K \log n$.
Thus we put into $A$, on the one hand, the points thrown away
during this process, and, on the other hand, the intervals for
which the return time does not reach $[3^p,3^{p+1}-1]$ before $K
\log n$ iterations. On all other intervals we define $\tau$ as the
first occasion when the return time belongs to $[3^p,3^{p+1}-1]$.
\begin{prop}
\label{prop_temps_arret_standard} The standard stopping time
$(A,\tau)$ defined this way is indeed a stopping time if $n$ is
large enough. Furthermore, $\Leb(A)/\Leb(\boC) \leq n^{-1/5}$.
\end{prop}
\begin{proof}
The only non-trivial condition to be verified is
$\Leb(A)/\Leb(\boC) \leq n^{-1/5}$.

Let us first estimate the measure of points thrown away during the
refinement process. We will denote this set by $A_0(\subset
A\subset \boC)$.

No matter which phase of the iteration we consider, the return
time is $\geq n^{1/4}$, thus, according to \eqref{mesure_jete},
the points thrown away occupy at most a $C n^{-1/4}$ proportion of
the considered interval. Hence, by bounded distortion, the
proportion of $A_0$ in $\boC$ is at most $C n^{-1/4} K \log n \leq
n^{-1/5}$ for $n$ large enough.

It remains to be shown that the overall measure of the intervals
that do not reach $[3^p,3^{p+1}-1]$ before $K\log n$ iterations is
small. We have $\boC=A_0 \cup \bigcup \boC_i$, where each $\boC_i$
is of the form $\boC_{n_0,\dots,n_k}$ for some admissible sequence
$n_0,\dots,n_k$, with $k\leq K \log n$, and $n_k <3^{p+1}$
whenever $k< \lfloor K \log n\rfloor$. Thus it is enough to
estimate the measure of $\boC_i$-s with $\tau_{\boC_i}=k= \lfloor
K \log n\rfloor$. Let $\boC'$ be one of our \emph{standard} curves
of return time $n$. We apply the same construction to $\boC'$, and
get a similar decomposition $\boC'=A'_0 \cup \bigcup \boC'_i$.
Furthermore, by Lemma~\ref{comparables}, $\frac{\Leb \boC_i}{\Leb
\boC'_i} \leq \Co$.

Recall that the standard curves of return time $n$ foliate the
major part $M'_n$ of the stripe $M_n$ (where
$\mu(M'_n)/\mu(M_n)=1+O(1/n)$). For a fixed $i=(n_0,\dots,n_k)$
consider $B_i$ the subset of the stripe $M_n$ that corresponds to
the union of such $\boC'_i$-s for all the standard curves of
return time $n$. As the density of $\mu$ on $M_n$ is bounded away
from $0$, we get $\frac{\Leb(\boC_i)}{\Leb(\boC)} \leq C \frac{
\mu(B_i)}{\mu(M'_n)}$. Fix $B$ as the union of all $B_i$-s with
$\tau_i=\lfloor K \log n\rfloor$. When pulled back to the Young
tower, the preimages of the points of $B$ are all at height at
least $K\log n$. As $\pi_X^*(\mu_{\bar \Delta})=\mu$, we get
$\mu(B)\leq C \rho^{K\log n}=O(1/n^4)$ by our choice of $K$. As
$\mu(M'_n)\sim C/n^3$, we may put all these estimates together to
conclude that
  \begin{equation*}
  \frac{ \sum_{\tau_i=\lfloor K \log n \rfloor} \Leb(\boC_i)}
  {\Leb(\boC)} =O(1/n).
  \end{equation*}
This completes the proof of the proposition.
\end{proof}

In the next proposition we consider standard curves $\boC$ and use
the notation $(A_\boC,\tau_\boC)$ for their standard stopping
times. We define a subset of the phase space, a suitable union of
subcurves of standard curves, as $Y=\bigcup_\boC (\boC\moins
A_\boC)$. We also consider the Birkhoff sum of $\gam$ with respect
to $T^{-1}$ up to standard stopping time, i.e., we fix
$\Gam(x)=\sum_{k=1}^{\tau_\boC(x)-1} \gam(T^{-k}x)$ for $x\in Y$.

\begin{prop}
\label{prop_reduit_a_X_1}
We have
  \begin{equation*}
  \int_{\bar\Delta} \bar\gam  (e^{it\bar\Gam}-1)
  = \int_{Y}
  \gam (e^{it\Gam}-1) +
  O(t).
  \end{equation*}
\end{prop}
This proposition plays a central role as it allows us to
investigate, instead of $\int_{\bar\Delta} \bar\gam
(e^{it\bar\Gam}-1)$ (a quantity that depends \emph{a priori} on
the choice of the Young tower), an expression which is much easier
to handle, as it is completely explicitly given in terms of the
phase space geometry.

\begin{proof}
Let us show first that
  \begin{equation}
  \label{premiere_etape}
  \int_{\bar \Delta} \bar\gam
  (e^{it\bar\Gam}-1) =
  \int_{\pi_X^{-1}(Y)} \bar\gam
  (e^{it\bar\Gam}-1) +O(t).
  \end{equation}
Consider $A=X\moins Y$. The set $A$ consists of two parts. It
contains, on the one hand, the points that are not covered by
standard curves and, on the other hand, those contained in
$A_\boC$ for some standard curve $\boC$. These two sets will be
referred to as $A_1$ and $A_2$, respectively.

We cover the set $A_1\cap \{\phi_-=n\}$ by two further sets, the
first one containing points that slide along a semi-circle (of
return time $n$, this is of measure $O(1/n^4)$), and secondly the
part of $M_n$ not covered by standard curves, this later having
measure $O(1/n^4)$ as well. Altogether we have $\mu(A_1 \cap
\{\phi_-=n\})=O(1/n^4)$.

According to Proposition~\ref{prop_temps_arret_standard}, we have
$\Leb(A_\boC)/\Leb(\boC)\leq n^{-1/5}$ whenever $n$, the return
time for $\boC$, is large enough. Integrating on the relevant
standard curves we obtain $\mu(A_2\cap
\{\phi_-=n\})=O(1/n^{3+1/5})$.

Altogether we have
  \begin{equation}
  \label{majore_mesure_A}
  \mu(A \cap \{\phi_-=n\})=O(1/n^{3+1/5}).
  \end{equation}
For any $1/p+1/q=1$ we have
  \begin{equation*}
  \left| \int_{\pi_X^{-1}(A)} \bar\gam
  (e^{it\bar\Gam}-1) \right|
  \leq \int 1_{\pi_X^{-1}(A)} |\bar\gam|  t|\bar \Gam|
  \leq |t| \left( \int (1_{\pi_X^{-1}(A)} |\bar\gam|)^p\right)^{1/p}
  \left( \int |\bar\Gam|^q \right)^{1/q}.
  \end{equation*}
Recall from Lemma~\ref{integrable} that the function $\bar\Gam$
belongs to $L^q$ for any $q<2$, while \eqref{majore_mesure_A}
implies that $\int (1_{\pi_X^{-1}(A)} |\bar\gam|)^p$, being equal
to $\int_X 1_A |\gam|^p$, is finite for $p<2+1/5$. We can thus
take $p=2+1/10$ and $q=(1-1/p)^{-1}$, to obtain
\eqref{premiere_etape}.

Now, to complete the proof, we need to show that
  \begin{equation}
  \label{deuxieme_etape}
  \int_{\pi_X^{-1}(Y)} \bar\gam
  (e^{it\bar\Gam}-1)
  =\int_{\pi_X^{-1}(Y)} \bar\gam
  (e^{it\Gam\circ
  \pi_X}-1)+O(t).
  \end{equation}

Consider $\boC_i$, a connected component of $\boC\moins A_\boC$,
where $\boC$ is a standard curve of return time $n$. Then the
stopping time on $\boC_i$ is an integer $\tau_i<K\log n$ such that
$\boD_i=T^{-\tau_i}(\boC_i)$ is a good curve, with return time in
the interval $[n^{1/4}/3,3 n^{1/4}]$.
\begin{lem}
\label{lemme_chernov} There exists a constant $\Cb$ such that, for
any large enough integer $n$, given any good curve $\boD$ of
return time $\in [n^{1/4}/3, 3 n^{1/4}]$, the points for which the
return time increases above $n^{1/2}$ within $K\log n$ iterations
of $T^{-1}$ occupy relative measure less than $\Cb n^{-1/4}$ in
$\boD$.
\end{lem}

\begin{proof}
The map $T^{-1}$ satisfies Chernov's axioms, by
\cite{markarian:slow}. Consequently, we can use \cite[Theorem
3.1]{chernov:decay}, with
$\delta=Z[\boD,\boD,0]^{-1/\sigma}/n^{1/\sigma}$. This theorem is
in fact stated for LUMs, but its proof can be straightforwardly
adapted to deal with manifolds close to the unstable direction.

We obtain a decreasing sequence $W_0^1 \supset W_1^1 \supset \dots
\supset W^1_{\lfloor K \log n \rfloor }$ of subsets of $\boD$ such
that, if we denote by $\Sing$ the set of singularities of
$T^{-1}$,
  \begin{equation}
  \label{eq1chernov}
  \forall c>0,
  \forall 0\leq p \leq K\log n,\ \Leb\{x\in W_p^1 \tq
  \dist(T^{-p}x,\Sing)\leq c n^{-1}\} \leq Cc n^{-1}
  \end{equation}
(by Equation (3.3) in \cite{chernov:decay}), and
  \begin{equation}
  \label{eq2chernov}
  \forall 0 \leq p \leq K\log n,\ \Leb(W_p^1 \moins W_{p+1}^1) \leq
\frac{C}{n} \Leb(\boD)
  \end{equation}
(By (iv), (3.5) in \cite{chernov:decay} and our choice of
$\delta$).

Note that the results of \cite{chernov:decay} imply that
\eqref{eq1chernov} holds for the distance measured in the
$p$-metric. However, we are in a region of $X$ where $\cos \theta$
is bounded away from $0$, and the stable and unstable cones are
bounded away from the vertical direction by Proposition
\ref{prop_hyperbolic}. Hence, it is equivalent to have
\eqref{eq1chernov} for the $p$-distance or for the usual distance.

If $T^{-p}(x)$ has a return time $\geq n^{1/2}$, then $T^{-p}x$ is
at a distance at most $C n^{-1}$ of $\Sing$. Hence, the point $x$
belongs to one of the sets whose measure is bounded in
\eqref{eq1chernov} and \eqref{eq2chernov}. This gives a measure at
most $C\log n\, n^{-1}$. Since $\Leb(\boD)\geq C n^{-1/2}$, this
proves the lemma.
\end{proof}

This lemma applies to $\boD_i$. Let us write $\boC_i=\boC_i^1 \cup
\boC_i^2$, where $\boC_i^2$ corresponds to points which go to
$\boD_i$, and then reach a return time $>n^{1/2}$ in a time
shorter than $K\log n$. It satisfies
$\Leb(\boC_i^2)/\Leb(\boC_i)\leq C n^{-1/4}$ by
Lemma~\ref{lemme_chernov}.

Let $Y_1= \bigcup \boC_i^1$ and $Y_2=\bigcup \boC_i^2$. Since
$\mu( Y_2 \cap \{\phi_-=n\}) = O(1/n^{3+1/4})$, the proof of
\eqref{premiere_etape} applies and gives
  \begin{equation*}
  \int_{\pi_X^{-1}(Y_2)} \bar\gam
  (e^{it\bar\Gam}-1)=O(t);\quad
  \int_{\pi_X^{-1}(Y_2)} \bar\gam
  (e^{it\Gam\circ
  \pi_X}-1)=O(t).
  \end{equation*}

\begin{rmq}
\label{H_in_lq}
Note that $\Gam\circ \pi_X$ belongs to $L^q$ for any $q<2$ as it
is smaller than a function to which Lemma~\ref{integrable}
applies.
\end{rmq}

Hence, it is sufficient to prove \eqref{deuxieme_etape} on
$\pi_X^{-1}(Y_1)$. Let us write $\pi_X^{-1}(Y_1)=Z_1\cup Z_2$
where
  \begin{equation*}
  Z_1=\{ x\in \pi_X^{-1}(Y_1) \tq \omega(x) < K \log( \phi_- (\pi_X x))\}
  \end{equation*}
and $Z_2=\pi_X^{-1}(Y_1) \moins Z_1$. For $n>0$,
  \begin{equation*}
  \mu_{\bar \Delta} \{ x\in Z_2 \tq \phi_-(\pi_X x)=n \}
  \leq \mu_{\bar \Delta} \{ x\in \bar \Delta, \omega(x)\geq K \log n\}
  =O(1/n^4).
  \end{equation*}
Hence, we get once again $\int_{Z_2} \bar\gam
(e^{it\bar\Gam}-1)=O(t)$ and $\int_{Z_2} \bar\gam (e^{it\Gam\circ
\pi_X}-1)=O(t)$.

On $Z_1\cap \{ \phi_- \circ \pi_X=n\}$, the functions $\bar\Gam$
and $\Gam \circ \pi_X$ differ by at most $\norm{f_0}_\infty K\log
n\, n^{1/2}$ (corresponding to at most $K\log n$ iterations with a
return time $<n^{1/2}$). Hence,
  \begin{multline*}
  \left| \int_{Z_1} \bar\gam
  (e^{it\bar\Gam}-1)
  - \bar\gam  (e^{it\Gam\circ \pi_X}-1) \right|
  \\
  \leq |t| \int_{Z_1} |\bar \gam|  |\bar\Gam-\Gam \circ \pi_X|
  \leq C|t| \sum_n \mu\{\phi_-=n\} n  \log n\, n^{1/2}
  \leq C|t|
  \end{multline*}
since $\mu\{\phi_-=n\}=O(1/n^3)$. This proves
\eqref{deuxieme_etape}, and concludes the proof of Proposition
\ref{prop_reduit_a_X_1}.
\end{proof}

\subsection{An upper bound on $\Gam$}

The aim of this subsection is to estimate the average of the
function $\Gam$ on a good curve ${\boC}$ of return time $n$. We
obtain the following upper bound:

\begin{prop}
\label{majore} Let $\rhp\in [1,2)$. Consider a good curve ${\boC}$
of return time $n_0$, and a stopping time $(A,\tau)$ on ${\boC}$.
Then
  \begin{equation*}
  \frac{\int_{{\boC} \moins A} \sum_{k=0}^{\tau(x)-1}
  |\gam(T^{-k}x)|^\rhp}{\Leb({\boC}\moins A)}
  \leq \Cd(\rhp) n_0^\rhp,
  \end{equation*}
where the constant $\Cd(\rhp)$ depends only on $\rhp$.
\end{prop}

\newcommand{\gamphi}{\phi_-}
\newcommand{\gamphib}{\phi'}

Let us fix some notation first. There is an integer $p_0$ such
that the return time $n_0$ for our good curve $\boC$ belongs to
$[3^{p_0}, 3^{p_0+1}-1]$. By the definition of stopping times,
there exists another integer $p_1<p_0$ such that, for any $x\in
\boC \moins A$, $\phi_-(T^{-\tau(x)}(x)) \in
[3^{p_1},3^{p_1+1}-1]$. Now consider an intermediate $p$, $p_1<
p\le p_0$. In the course of the proof first we investigate, in a
series of lemmas, what happens while the return time descends from
$[3^p,3^{p+1}-1]$ to $[3^{p-1},3^p-1]$. Then we sum up for $p_1<
p\le p_0$. In the first part of the proof the value of $p$ is
fixed and $n\approx 3^p$, while in the second part $p$ varies from
$p_1$ to $p_0$. The value of $\rhp\in [1,2)$ is fixed throughout
the subsection.

According to this plan, let us fix $p\in\N$ large enough. Given
$x\in X$, we define $\tau_{p}(x)$ as the first time $k\geq 1$ for
which $\phi_-(T^{-k}x)<3^p$, and $\Phi_p(x)
=\sum_{k=0}^{\tau_{p}(x)-1} |\phi_-(T^{-k}x)|^{\rhp}$. Since
$|\gam|\leq C \phi_-$, it is sufficient to prove Proposition
\ref{majore} for $h=\phi_-$ to conclude.

Define $R\subset X$ as the union of all standard curves with
return time from the interval $[3^p/2, 3^p-1]$.
\begin{lem}
\label{majore_sur_R}
There exists a constant $\Cz$ such that
  \begin{equation*}
  \int_R \Phi_p \leq \Cz \mu(R) 3^{ps}.
  \end{equation*}
\end{lem}
\begin{proof}
Let $R_1=\{x\in R \tq \phi_-(T^{-1}x)<3^p\}$ and $R_2=\{x\in R \tq
\phi_-(T^{-1}x) \geq 3^p\}$. On $R_1$ we have
$\Phi_p(x)=|\gamphi(x)|^\rhp$, thus
  \begin{equation*}
  \int_{R_1} \Phi_p \leq C \mu(R_1) 3^{ps}.
  \end{equation*}

Let us define $\gamphib(x) =\gamphi(x)$ for $x$ with $\phi_-(x)
\geq 3^{p-1}$ and $\gamphib(x)=0$ otherwise. Note that $\Phi_p(x)
=\sum_{k=0}^{\tau_{p}(x)-1} |\gamphib(T^{-k}x)|^{\rhp}$ for $x\in
R_2$.

Consider $Z\subset X$, $Z:= \{3^{p-1}-\Cmm \leq \phi_- <3^p\}$,
and define $\tau_Z:Z\to \N$ as the first return time to $Z$. By
Kac's formula,
  \begin{equation*}
  \int_Z \sum_{k=0}^{\tau_Z(x)-1} |\gamphib(T^{-k}x)|^\rhp =
  \int_X |\gamphib|^s \leq C \sum_{k\geq 3^{p-1} } \mu(\phi_-=k) |k|^\rhp
  \leq C \sum_{k \geq 3^{p-1}} \frac{1}{k^3} k^\rhp \leq C
  \frac{3^{p\rhp}}{3^{2p}}.
  \end{equation*}
Now  $R_2\subset Z$ and for $x\in R_2$ we have
$\tau_Z(x)=\tau_p(x)$. Thus
  \begin{equation*}
  \int_{R_2} \Phi_p = \int_{R_2}
  \sum_{k=0}^{\tau_Z(x)-1} |\gamphib(T^{-k}x)|^\rhp
  \leq \int_Z \sum_{k=0}^{\tau_Z(x)-1} |\gamphib(T^{-k}x)|^\rhp.
  \end{equation*}
By Remark~\ref{transition_proba}, $\frac{1}{3^{2p}}=O(\mu(R_2))$.
This completes the proof.
\end{proof}

If $\boC$ is a good curve of return time $n \in [3^p,3^{p+2}-1]$,
$\tau_p$ defines a stopping time on $\boC$, with the corresponding
thrown-away set that we denote by $A_p$. To see that it is indeed
a stopping time we only need to show that $\Leb(A_p)\leq
\Leb(\boC)/2$. Now consider the standard stopping time
$\tau_{\boC}$ with its thrown away set $A_\boC$. Then $A_p\subset
A_\boC$ while $\Leb(A_\boC) \leq n^{-1/5}\Leb(\boC)$ by
Proposition~\ref{prop_temps_arret_standard}, which gives the
claim.

The first step in the proof of Proposition~\ref{majore} is the
estimate
  \begin{equation}
  \label{eq_one_step}
  \frac{ \int_{\boC \moins A_p} \Phi_p}{\Leb(\boC \moins A_p)} \leq C
  3^{p\rhp}.
  \end{equation}
for a good curve $\boC$ with return time $n\in [3^p,3^{p+1}-1]$.
To show this, we will relate the average of $\Phi_p$ on $\boC$ to
its average on $R$.

Consider $B=\bigcup (\boC \moins A_p)$, where the union is taken
over all standard curves of return time from the interval $[3^p,
3^{p+2}-1]$.

\begin{lem}
\label{ineg_un_sens} There is a constant $\Cy$ such that, for any
good curve $\boC$ of return time $n\in [3^p,3^{p+1}-1]$,
  \begin{equation}
  \label{alshdfljhasdf}
  \frac{ \int_{\boC \moins A_{p}} \Phi_p}{\Leb(\boC \moins A_p)} \leq
  \Cy \frac{\int_B \Phi_p}{\mu(B)} + \Cy 3^{p\rhp}.
  \end{equation}
\end{lem}
\begin{proof}
Let $U=\{x\in \boC \tq \phi_-(T^{-1}x) \geq 3^p\}$. On $\boC$, we
have $\Phi_p(x) = |\gamphi(x)|^s + 1_U(x) \Phi_p(T^{-1}x)$. To
prove \eqref{alshdfljhasdf}, it is enough to show
  \begin{equation*}
  \int_{ U\cap (\boC \moins A_p)} \Phi_p\circ T^{-1} \leq C
  \frac{\int_B \Phi_p}{\mu(B)} \Leb(\boC \moins A_p).
  \end{equation*}
By bounded distortion, this can be further reduced to
  \begin{equation}
  \label{asjhdfkljashdf}
  \frac{\int_{ T^{-1}(\boC \moins A_p) \cap \{\phi_- \geq 3^p\}} \Phi_p}
  {\Leb(  T^{-1}(\boC \moins A_p) \cap \{\phi_- \geq 3^p\})}
  \leq C \frac{\int_B \Phi_p}{\mu(B)}.
  \end{equation}

\begin{comment}
Here we need distortion bound for $T^{-1}$ on all $\boC$. This is
OK, though it would not be true for a (much) higher iterate
$T^{-k}$.
\end{comment}

Let $q$ be the maximal possible return time the points of
$T^{-1}(\boC \moins A_p)$ have. It satisfies $3^{p+2} > q \geq
3^{p+1}-\Cn$. By Lemma~\ref{comparables}
  \begin{equation*}
  \frac{\int_{ T^{-1}(\boC \moins A_p) \cap \{\phi_- \geq 3^p\}} \Phi_p}
  {\Leb(  T^{-1}(\boC \moins A_p) \cap \{\phi_- \geq 3^p\})}
  \leq C \frac{ \int_{B\cap \{3^p \leq
  \phi_- \leq q\}}\Phi_p}{ \mu(B\cap \{3^p \leq
  \phi_- \leq q\})}.
  \end{equation*}

As $q\geq 3^{p+1}-\Cn$, by Remark~\ref{transition_proba} $\mu(B)
\leq C \mu(B\cap \{3^p \leq \phi_- \leq q\})$. This implies
\eqref{asjhdfkljashdf} and completes the proof.
\end{proof}

Now $B$ is not exactly $R$, we need to ``widen up'' the estimate
of Lemma~\ref{majore_sur_R} from $R$ to $B$ to obtain
\eqref{eq_one_step}.

Let $B_1 = B\cap\{3^p\leq \phi_- < 3^{p+1}/2\}$, $B_2=B
\cap\{3^{p+1}/2 \leq \phi_- < 3^{p+1}\}$, $B_3=B\cap\{3^{p+1}\leq
\phi_- < 3^{p+2}/2\}$ and $B_4 =B\cap \{3^{p+2}/2 \leq \phi_-
<3^{p+2}\}$.
\begin{lem}
\label{autre_ineg} There exists a constant $\Cx$ such that, for
any good curve $\boC$ of return time $n\in [3^p,3^{p+1}/2)$,
  \begin{equation*}
  \frac{ \int_{B_1} \Phi_p}{\mu(B_1)} \leq \Cx \frac{\int_{\boC \moins
  A_p} \Phi_p}{\Leb(\boC \moins A_p)}.
  \end{equation*}
\end{lem}
\begin{proof}
The curve $T^{-1}(\boC)$ crosses all stripes of return time
between $3^p$ and $3^{p+1}/2$. This allows us to apply the
argument of Lemma~\ref{ineg_un_sens} with reversed inequalities.
\end{proof}

\begin{lem}
\label{majore_sur_B1} There exists a constant $\Cw$ such that
  \begin{equation*}
  \frac{\int_{B_1} \Phi_p}{\mu(B_1)} \leq \Cw 3^{p\rhp}.
  \end{equation*}
\end{lem}
\begin{proof}
Let $\boC$ be a standard curve of return time $n\in
[3^p/2,3^p-1]$. For $i\in [3^p, 3^{p+1}/2)$, put $\boC_i =\{x\in
\boC, \phi_-(T^{-1}x)=i\}$ and let  $\boD_i$ be its image by
$T^{-1}$. This is a good curve of return time $i$ and, by bounded
distortion,
  \begin{equation*}
  \frac{\int_{\boD_i} \Phi_p}{\Leb(\boD_i)} \leq C
  \frac{\int_{\boC_i}\Phi_p}{\Leb( \boC_i)}.
  \end{equation*}
Furthermore, applying Lemma~\ref{autre_ineg} to $\boD_i$,we get
  \begin{equation*}
  \Leb(\boC_i) \frac{ \int_{B_1} \Phi_p}{\mu(B_1)} \leq C
  \int_{\boC_i}\Phi_p.
  \end{equation*}
As by Remark~\ref{transition_proba} the good curves $\boC_i$
occupy a fixed proportion of $\boC$, we may sum up
  \begin{equation*}
  \Leb(\boC) \frac{ \int_{B_1} \Phi_p}{\mu(B_1)} \leq C \int_{\boC}
  \Phi_p.
  \end{equation*}
Integrating over all standard curves of return time $\in
[3^p/2,3^p-1]$, we obtain
  \begin{equation*}
  \mu(R) \frac{ \int_{B_1} \Phi_p}{\mu(B_1)} \leq C \int_R \Phi_p.
  \end{equation*}
We may conclude by Lemma~\ref{majore_sur_R}.
\end{proof}

\begin{lem}
\label{majore_sur_B2} There is a constant $\Cv$ such that for any
$l=2,3,4$,
  \begin{equation*}
  \frac{\int_{B_l} \Phi_p}{\mu(B_l)} \leq \Cv 3^{p\rhp}.
  \end{equation*}
\end{lem}
\begin{proof}
As the three cases are essentially identical we give the argument
only for one of them, for $l=3$, say. The proof is analogous to
that of the previous lemma, we only need to apply a bit more
iterations. Let $\boC$ be a standard curve with return time from
$[3^p/2,3^p-1]$. Given $i\in [3^p,3^{p+1}/2)$, let $\boC_i$ be the
set of points in $\boC$ the images of which have return time $i$.
For $j\in [3^{p+1}/2, 3^{p+1})$, let $\boC_{ij}$ be the set of
points in $\boC_i$ the $T^{-2}$-images of which have return time
$j$. Finally, for $k\in [3^{p+1},3^{p+2}/2)$, we define
$\boC_{ijk}$ analogously.

By Remark~\ref{transition_proba}, at each step we keep a fixed
proportion of the previous set. Thus, there exists a constant $C$
such that
  \begin{equation*}
  \Leb(\boC) \leq C \sum_{i,j,k} \Leb(\boC_{ijk}).
  \end{equation*}
Following the lines of the proof of Lemma~\ref{autre_ineg} we may
show that given any good curve $\boD$ of return time from the
interval $[3^{p+1},3^{p+2}/2)$, we have $\frac{\int_{B_3}
\Phi_p}{\mu(B_3)} \leq C \frac{\int_\boD \Phi_p}{\Leb(\boD)}$.
This applies, in particular, to $\boD =T^{-3}(\boC_{ijk})$ and
gives
  \begin{equation*}
  \frac{\int_{B_3} \Phi_p}{\mu(B_3)}
  \leq C \frac{\int_{T^{-3} \boC_{ijk}}
  \Phi_p}{\Leb(T^{-3}\boC_{ijk})}
  \leq C \frac{ \int_{\boC_{ijk}} \Phi_p}{\Leb(\boC_{ijk})},
  \end{equation*}
by bounded distortion. We may apply Lemma~\ref{majore_sur_R}, just
as we did in the proof of Lemma~\ref{majore_sur_B1}, to get the
desired conclusion.
\end{proof}

Lemmas~\ref{majore_sur_B1}, \ref{majore_sur_B2} and
\ref{ineg_un_sens} altogether imply the bound \eqref{eq_one_step}
for any good curve of return time $n\in [3^p,3^{p+1}-1]$. We apply
this bound in the second (much easier) step of the proof of
Proposition~\ref{majore}.

\begin{proof}[Proof of Proposition~\ref{majore}]
Recall the notations from the beginning of the subsection: $\boC$
is a good curve of return time $n_0\in [3^{p_0},3^{p_0+1}-1]$, for
some large $p_0$, and the stopping time $\tau$ is related to
another integer $p_1$ ($p_0>p_1$): $\phi_-(T^{-\tau(x)}(x)) \in
[3^{p_1},3^{p_1+1}-1]$ for all $x\in \boC \moins A$.

To simplify notation in this proof we define $\tau_{p_0+1}(x)=0$
for $x\in \boC \moins A$. For $x\in \boC\moins A$ we have
  \begin{equation*}
  \sum_{k=0}^{\tau(x)-1}|\gamphi(T^{-k}x)|^\rhp = \sum_{p=p_1+1}^{p_0}
  \Phi_p(T^{-\tau_{p+1}(x)}x).
  \end{equation*}

Let $p_1+1\leq p\leq p_0$ and $x\in {\boC}\moins A$. Then there is
a subcurve ${\boC}_i \subset \boC$ that contains $x$ and for which
$T^{-\tau_{p+1}(x)}({\boC}_i)$, to be denoted by $\boD_i$, is a
good curve of return time from $[3^p,3^{p+1}-1]$. By bounded
distortion
  \begin{equation*}
  \frac{\int_{{\boC}_i \moins A} \Phi_{p}(T^{-\tau_{p+1}}y)}
  {\Leb({\boC}_i )}
  \leq C \frac{\int_{T^{-\tau_{p+1}}(\boC_i\moins A)} \Phi_p}
  {\Leb(\boD_i)}
  \leq C \frac{\int_{\boD_i \moins A_p} \Phi_p}{\Leb(\boD_i)}.
  \end{equation*}
Now according to \eqref{eq_one_step} this final quantity is
bounded from above by $C 3^{ps}$. Summing up for all intervals
$\boC_i$ we obtain
  \begin{equation*}
  \int_{\boC\moins A}\Phi_p(T^{-\tau_{p+1}(y)}y) \leq C 3^{ps}
  \Leb(\boC).
  \end{equation*}
Summation on $p$ from $p_1+1$ to $p_0$ implies the statement.
\end{proof}

\begin{cor}
\label{cor_sans_exp} We have
  \begin{equation*}
  \int_{Y} \gam  (e^{it\Gam}-1)
  = it \int_{Y} \gam 1_{\phi_- \leq 1/|t|}\Gam
  +o(t \log(1/|t|)).
  \end{equation*}
\end{cor}
\begin{proof}
We have
  \begin{equation*}
  \left| \int_{Y} \gam 1_{\phi_- \leq 1/|t|} (e^{it\Gam}-1-it\Gam)\right|
  \leq C \int_{Y } |\gam| 1_{\phi_- \leq 1/|t|} |t|^{3/2} |\Gam|^{3/2}.
  \end{equation*}
We may estimate $|\Gam(x)|^{3/2}$ as
  \begin{equation*}
  |\Gam(x)|^{3/2}\leq \left(\sum_1^{\tau(x)-1} |\gam(T^{-k}x)|\right)^{3/2}
  \leq \tau(x)^{1/2} \sum_0^{\tau(x)-1} |\gam(T^{-k}x)|^{3/2}.
  \end{equation*}
Now put $\Phi(x)=\sum_0^{\tau(x)-1}|\gam(T^{-k}x)|^{3/2}$. Then
for $x\in Y$ of return time $n$ we get $|\Gam(x)|^{3/2} \leq
(K\log n)^{1/2} \Phi(x)$, as the standard stopping time satisfies
$\tau(x)\leq K \log n$.

By Proposition~\ref{majore} the average of the function $\Phi$ on
$Y \cap \{\phi_-=n\}$ is less than $c n^{3/2}$. Putting these
estimates together
  \begin{align*}
  \int_{Y} |\gam| 1_{\phi_- \leq 1/|t|} |t|^{3/2} |\Gam|^{3/2} &
  \leq C |t|^{3/2} \sum_{n=1}^{1/|t|} \mu(\phi_-=n)  \sqrt{\log n}\, n\,
  n^{3/2}
  \\&
  \leq C |t|^{3/2}  \sqrt{ \log(1/|t|)} |t|^{-1/2}
  =o( |t| \log(1/|t|)),
  \end{align*}
while
  \begin{equation*}
  \left| \int_Y \gam 1_{\phi_->1/|t|} (e^{it\Gam}-1)\right|
  \leq C \int \phi_- 1_{\phi_->1/|t|}\leq C \sum_{n>1/|t|}
  \mu(\phi_-=n) n= O(t)
  \end{equation*}
as $\mu(\phi_-=n)=O(1/n^3)$.
\end{proof}

\subsection{Exact asymptotics for $\Gam$}

Recall the value of $I$ from \eqref{donne_I}, and the fact that on
$Y\cap \{\phi_-=n\}$ the function $h$ is equivalent to $nI$.

\begin{lem}
\label{lemme_estimee_precise} Let $y=\frac{1}{1-\frac{3}{4}\log
3}$. For any $\epsilon>0$ there exists $N_0\in\N$ such that, for
all $n\geq N_0$, for all good curve $\boC$ with return time $n$,
  \begin{equation*}
  \left|\frac{\int_{\boC \moins A_\boC} H}{\Leb(\boC\moins A_\boC)}-n(y-1)I
  \right| \leq \epsilon n.
  \end{equation*}
\end{lem}
\begin{proof}

Recall the asymptotic expressions for the transition probabilities
from Remark~\ref{transition_proba}. These allow us to regard the
map $T^{-1}$ as a Markov chain. Then the statement of the lemma
can be guessed by the expectation value with respect to the
invariant distribution of this chain.

The rigorous proof is inductive. Note that first we fix
$\epsilon>0$, that will correspond to the required precision in
the asymptotics, and then we may choose $n$ arbitrarily large. Let
$L\in \N$ be an integer for which $(9/10)^L \leq \epsilon$. This
integer $L$ is the number of inductive steps needed to obtain
$\epsilon$-precision. More precisely, if $n_0,\ldots,n_L$ is an
admissible sequence (here $n_0=n$), then $n_L$ is typically much
smaller than $n_0$. The Birkhoff sum of $\gam$ for the times
between $n_L$ and the stopping time can be estimated by the upper
bound coming from Lemma~\ref{majore}, which roughly means that we
only need to take care of the sum for the first $L$ steps. This
estimate will be the starting point of our induction. Then we
place our standard curve ``high enough'' (i.e., choose $n$ large
enough) to ensure that the transition probabilities of
Remark~\ref{transition_proba} are accurate with very good
precision. These transition probabilities are responsible for the
appearance of $y$ as we decrease the length of the admissible
sequence $n_0,\ldots,n_i$ from $i=L$ to $i=0$ in the induction.

Let $\boC$ be a standard curve of return time $n$ with the
standard stopping time $(A_\boC,\tau_\boC)$ on it. If
$n_0,\dots,n_i$ is admissible with $n_0=n$ and $i\leq L$, the set
${\boC}_{n_0,\dots,n_i}$ is not empty, and we may consider
$\boC'_{n_0,\ldots,n_i}=\boC_{n_0,\ldots,n_i} \cap (\boC\moins
A_\boC)$. For $\boD=T^{-i}(\boC_{n_0,\ldots,n_i})$, define $A=\boD
\moins T^{-i}(\boC'_{n_0,\ldots,n_i})$ and
$\tau(T^{-i}x)=\tau_\boC(x)-i$. Then, for large enough $n$,
$(A,\tau)$ is a stopping time on $\boD$. To see this we note that
$\Leb(A)\leq \Leb(\boD)/2$ as the number of iterations is bounded
from above by $L$ while $\Leb(A_\boC)/\Leb(\boC) \to 0$ as
$n\to+\infty$.

By increasing $n$ if necessary, we may assume that for any
$p>n/3^L$, and for any $x$ with $\phi_-(x)=p$ we have $|h(x)-p I|
\leq p/(L 3^L)$.

Thus for $x\in \boC'_{n_0,\dots,n_{L}}$ we have
  \begin{equation}
  \label{majore_H}
  | H(x)- (n_1+\dots+n_{L-1})I|
  \leq \sum_{k=1}^{L-1} |h(T^{-k}x)-n_k I| +
  \sum_{k=L}^{\tau_\boC(x)-1} |h(T^{-k}(x)|
  \end{equation}
where the first term satisfies
  \begin{equation*}
  \sum_{k=1}^{L-1} |h(T^{-k}x)-n_k I| \leq \sum_{k=1}^{L-1}
n_k/(L 3^L) \leq  n_L,
  \end{equation*}
as $n_k \leq 3^L n_L$. On the other hand if we integrate the
second term in \eqref{majore_H}, we may use the upper bound of
Lemma~\ref{majore}. We get, for some constant $\Ce$:
  \begin{equation}
  \label{initie_recurrence}
  \left|\frac{\int_{\boC'_{n_0,\dots,n_L}} \Gam- (n_1+\dots+n_{L-1})I}
  {\Leb(\boC'_{n_0,\dots,n_L})} \right| \leq \Ce n_L.
  \end{equation}

Choose $n$ large enough to ensure that (i) all the $\epsilon_p$
from Remark~\ref{transition_proba} are less than $\epsilon$
whenever $p>n/3^L$, and that (ii) the distortion of any
$T^{-i}\mid_{{\boC}_{n_0,\dots,n_i}}$, $i\le L$ is bounded from
above by $\epsilon$.

As $y\frac{3}{4}\log 3 -y+1=0$ we have, for $n$ large enough,
  \begin{equation}
  \label{caract_y}
  \left| y\sum_{p/3+\Cn}^{3p-\Cn} \frac{3}{8k}-y+1 \right| <\epsilon
  \end{equation}
whenever $p>n/3^L$.

To simplify notation we introduce $\alpha=\frac{9}{10}$ and
another positive number, $\beta> 2\log 3$ which is, however, not
too big so that $\frac{3}{8}\beta <\alpha$. By further increasing
$n$, if necessary, we may also assume that $\sum_{p/3}^{3p}
\frac{1}{k} \leq \beta$ whenever $p>n/3^L$.

Now, by induction on decreasing $i$ we show the following bound:
  \begin{equation}
  \label{recurrence_terrible}
  \left|\frac{ \int_{\boC'_{n_0,\dots,n_i}}
  \Gam-(n_1+\dots+n_{i-1})I}{n_i
  \Leb(\boC'_{n_0,\dots,n_i})}-yI\right|  \leq \alpha^{L-i} (\Ce+y|I|)+
  \Cf \sum_{k=i}^{L-1}
  \epsilon  \alpha^{k-i},
  \end{equation}
where $\Cf$ is some constant. Note that for $i=0$, when the sum
$n_1+\dots+n_{i-1}$ is to be interpreted as $-n_0$, this bound
implies the statement of Lemma~\ref{lemme_estimee_precise}. On the
other hand, the case $i=L$ is already established in
\eqref{initie_recurrence}. So let us assume
\eqref{recurrence_terrible} holds for $i$, and show it for $i-1$.
We have
  \begin{align*}
  \frac{ \int_{\boC'_{n_0,\dots,n_{i-1}}} \Gam-(n_1+\dots+n_{i-2})I}{n_{i-1}
  \Leb(\boC'_{n_0,\dots,n_{i-1}})} -yI
  \!\!\!\!\!\!\!\!\!\!\!\!\!\!\!\!\!\!\!\!\!\!\!\!\!\!\!\!\!\!
  \!\!\!\!\!\!\!\!\!\!\!\!\!\!\!\!\!\!\!\!\!\!\!\!\!\!\!\!\!\!
  &
  \\&
  =\frac{ \int_{\boC'_{n_0,\dots,n_{i-1}}} \Gam-(n_1+\dots+n_{i-1})I}{n_{i-1}
  \Leb(\boC'_{n_0,\dots,n_{i-1}})} -yI +I
  \\&
  =\sum_{n_i=n_{i-1}/3+\Cn}^{3n_{i-1}-\Cn}
  \frac{ \int_{\boC'_{n_0,\dots,n_{i}}} \Gam-(n_1+\dots+n_{i-1})I}{n_{i-1}
  \Leb(\boC'_{n_0,\dots,n_{i-1}})} -yI +I
  \\&
  =\sum_{n_i=n_{i-1}/3+\Cn}^{3n_{i-1}-\Cn} \left(
  \frac{ \int_{\boC'_{n_0,\dots,n_{i}}} \Gam-(n_1+\dots+n_{i-1})I}{n_{i}
  \Leb(\boC'_{n_0,\dots,n_{i}})} \frac{n_i \Leb(\boC'_{n_0,\dots,n_{i}})}{n_{i-1}
  \Leb(\boC'_{n_0,\dots,n_{i-1}})} - \frac{3}{8n_i} yI\right)
  \\ & \hphantom{=\ }
  +\left(y \sum_{n_i=n_{i-1}/3+\Cn}^{3n_{i-1}-\Cn} \frac{3}{8 n_i}
   -y +1\right)I.
  \end{align*}

The choice of a large enough $n$ ensures that even $n_i$ is large
enough so that \eqref{caract_y} applies:
  \begin{equation*}
  \left| y \sum_{n_i=n_{i-1}/3+\Cn}^{3n_{i-1}-\Cn} \frac{3}{8 n_i}
   -y +1\right| \leq \epsilon.
  \end{equation*}

Now we will use the transition probabilities
\eqref{proba_transition} on the curve
$T^{-(i-1)}(\boC_{n_0,\dots,n_{i-1}})$. We will also use that the
distortions of $T^{-(i-1)}$, when restricted to this curve, are
bounded from above by $\epsilon$. Note furthermore that
$\boC'_{n_0,\dots,n_{i-1}}$ occupies at least
$(1-\epsilon)$-proportion of $\boC_{n_0,\dots,n_{i-1}}$ if $n$ is
large enough (we may apply
Proposition~\ref{prop_temps_arret_standard}). The same holds for
$\boC'_{n_0,\dots,n_{i}}$ in $\boC_{n_0,\dots,n_i}$. These
observations allow us to obtain (note $n_{i-1}/n_i\le 3$):
  \begin{equation*}
  \left| \frac{\Leb(\boC_{n_0,\dots,n_i})}{\Leb(\boC_{n_0,\dots,n_{i-1}})} -
  \frac{\Leb(\boC'_{n_0,\dots,n_i})}{\Leb
  (\boC'_{n_0,\dots,n_{i-1}})} \right| \leq
  2\epsilon \frac{ \Leb(\boC_{n_0,\dots,n_i})}{\Leb(\boC_{n_0,\dots,n_{i-1}})}
  \leq \frac{C\epsilon}{n_i}.
  \end{equation*}
One more reference to Remark~\ref{transition_proba} and to the
fact that the distortions can be made smaller than $\epsilon$ if
$n$ is large enough implies
  \begin{equation*}
  \left|\frac{n_i \Leb(\boC'_{n_0,\dots,n_{i}})}{n_{i-1}
   \Leb(\boC'_{n_0,\dots,n_{i-1}})} - \frac{3}{8n_i} \right| \leq
  \frac{\Cg \epsilon}{n_i}.
  \end{equation*}
By the triangular inequality,
  \begin{align*}
  \Biggl|
  \frac{ \int_{\boC'_{n_0,\dots,n_{i}}} \Gam-(n_1+\dots+n_{i-1})I}{n_{i}
  \Leb(\boC'_{n_0,\dots,n_{i}})}\!\!&\  \frac{n_i
  \Leb(\boC'_{n_0,\dots,n_{i}})}{n_{i-1}
   \Leb(\boC'_{n_0,\dots,n_{i-1}})} - \frac{3}{8n_i} yI\Biggr|
  \\&
  \leq \left| \frac{ \int_{\boC'_{n_0,\dots,n_{i}}} \Gam-
  (n_1+\dots+n_{i-1})I}{n_{i}
  \Leb(\boC'_{n_0,\dots,n_{i}})} -yI \right| \frac{3}{8 n_i}
  \\& \ \
  +  \left| \frac{ \int_{\boC'_{n_0,\dots,n_{i}}} \Gam-
  (n_1+\dots+n_{i-1})I}{n_{i}
  \Leb(\boC'_{n_0,\dots,n_{i}})} \right| \left| \frac{n_i
  \Leb(\boC'_{n_0,\dots,n_{i}})}{n_{i-1}
   \Leb(\boC'_{n_0,\dots,n_{i-1}})} - \frac{3}{8n_i} \right|.
  \end{align*}
Let $B_i$ be the bound at step $i$ of the induction. Then the
first term is bounded from above by $\frac{3 B_i}{8 n_i}$, and the
second term is bounded from above by $\frac{(B_i+y|I|) \Cg
\epsilon}{n_i}$.

Recall the definitions of $\alpha$ and $\beta$, we have
$\sum_{p/3}^{3p} \frac{1}{k} \leq \beta$ and, if $\epsilon$ is
small enough, $(\frac{3}{8}+\Cg\epsilon)\beta<\alpha$.

Putting our estimates together we get
  \begin{align*}
  B_{i-1} & = \epsilon|I|+ \sum_{n_{i-1}/3+\Cn}^{3n_i-\Cn}
  \left[ \frac{3 B_i}{8 n_i}
  + \frac{(B_i+y|I|) \Cg \epsilon}{n_i}\right]
  \leq \epsilon|I|+ \left( \frac{3B_i}{8}+
  (B_i+y|I|) \Cg \epsilon \right)\beta
  \\&
  \leq (|I|+\Cg y|I|\beta)\epsilon + \alpha B_i,
  \end{align*}
Now if  \eqref{recurrence_terrible} holds for $i$ with $\Cf=
|I|+\Cg y |I| \beta$, it holds for $i-1$ with the same constants.

Taking $i=0$ we get
  \begin{equation*}
  \frac{1}{n} \, \left| \frac{\int_{\boC\moins A_\boC} \Gam}{\Leb({\boC\moins
  A_\boC})} -n (y-1)I \right|
  \leq C\alpha^L+C \epsilon
  \leq \Ch\epsilon
  \end{equation*}
by the choice of $L$. Note that the constant $\Ch$ depends only on
$I$, thus it can be ``swallowed'' by $\epsilon$. This completes
the proof of the lemma.
\end{proof}

\begin{prop}
\label{final_for_H}
We have
  \begin{equation*}
  \int_Y \gam 1_{\phi_-\leq 1/|t|}\Gam= \left(\frac{I^2
  (y-1)\ell^2}{\pi}+o(1)\right) \log(1/|t|).
  \end{equation*}
\end{prop}
\begin{proof}
First let us show that
  \begin{equation}
  \label{o_sans_phi1}
  \int_Y (\gam-\phi_- I) 1_{\phi_- \leq 1/|t|} \Gam=o(\log(1/|t|)).
  \end{equation}
Fix $\epsilon>0$. If $N$ is large enough we have $|\gam-\phi_-
I|\leq \epsilon \phi_-$ for $\phi_- \geq N$. Thus we get (note
that $\Gam$ is integrable, cf.\ Remark~\ref{H_in_lq})
  \begin{equation*}
  \left|\int_Y (\gam-\phi_- I) 1_{\phi_- \leq 1/|t|} \Gam\right|
  \leq O(1)+\sum_{N \leq n \leq 1/|t|} \epsilon n
  \int_{Y\cap \{\phi_-=n\}} |\Gam|.
  \end{equation*}
We may apply Proposition~\ref{majore} with $\rhp=1$ to show
$\int_{Y\cap \{\phi_-=n\}} |\Gam| \leq Cn \mu(\phi_-=n)=O(1/n^2)$.
Thus we get
  \begin{equation*}
  \left|\int_Y (\gam-\phi_- I) 1_{\phi_- \leq 1/|t|} \Gam\right| \leq
  O(1)+ C \epsilon \log(1/|t|) \leq C' \epsilon \log(1/|t|) .
  \end{equation*}
As the above inequality is true for any fixed $\epsilon>0$, we get
\eqref{o_sans_phi1}.

Now we estimate
  \begin{equation*}
  \int_Y \phi_- I 1_{\phi_- \leq 1/|t|} \Gam = \sum_{n=1}^{1/|t|}
  nI\int_{ Y \cap \{\phi_-=n\}} \Gam.
  \end{equation*}
By Lemma~\ref{lemme_estimee_precise} we have $\int_{ Y \cap
\{\phi_-=n\}} \Gam \sim (y-1)In \mu(\phi_-=n) \sim (y-1)In
\frac{\ell^2}{\pi n^3}$. Actually, the measure of the set
$\{\phi_-=n\}$ can be estimated by direct geometric arguments. Up
to negligible terms, it is equivalent to $\frac{\ell^2}{4\pi n^3}$
in all relevant zones of $X$ which are ``corners of
parallelograms''. As there are $4$ such relevant zones we obtain
the above formula.

Finally we get
  \begin{equation*}
   \int_Y \phi_- I 1_{\phi_- \leq 1/|t|} \Gam
  \sim \sum_{n=1}^{1/|t|} n I^2 (y-1) n \frac{\ell^2}{\pi n^3}
  \sim \frac{I^2 (y-1)\ell^2}{\pi} \log(1/|t|),
  \end{equation*}
which completes the proof.
\end{proof}

Proposition~\ref{estime_integrale_precise} follows from the
combination of Proposition~\ref{final_for_H},
Corollary~\ref{cor_sans_exp}, Proposition~\ref{prop_reduit_a_X_1}
and  Lemmas~\ref{lem_tower2}, \ref{lem_tower1}.

\section{Proof of the main theorems}

In this section, we prove Theorems~\ref{main_theorem} and
\ref{main_thm_2}. The main tool will be an abstract theorem
showing that, if an induced map satisfies a limit theorem, then
the original map satisfies the same limit theorem. Such a result
has been proved in the case of flows by \cite{melbourne_torok},
and extended to the discrete time case (and to non-polynomial
normalizations) in \cite{gouezel:skewproduct}. For the convenience
of the reader, we state here the result we will use.

If $Y$ is a subset of a probability space $(X,m)$, $T:X\to X$, and
$T_Y$ is the induced map on $Y$, we will write $S_n^Y
g=\sum_{k=0}^{n-1} g\circ T_Y^k$: this is the Birkhoff sum of $g$,
for the transformation $T_Y$. We will also write
$E_Y(g)=\frac{\int_Y g}{m[Y]}$. Finally, for $t\in \R$, $\lfloor t
\rfloor$ denotes the integer part of $t$.

\begin{thm}
\label{thm_probabiliste_general}
Let $T:X\to X$ be an ergodic endomorphism of a probability space
$(X,m)$, and $f:X\to \R$ an integrable function with vanishing
integral. Let $Y\subset X$ have positive measure. For $y\in Y$,
write $\phi(y)=\inf\{n>0 \tq T^n(y)\in Y\}$ and
$f_Y(y)=\sum_{k=0}^{\phi(y)-1} f(T^k y)$.

We assume the following properties:
\begin{enumerate}
\item
There exists a sequence $B_n\to +\infty$, with $\inf_{r \geq n}
\frac{B_r}{B_n}>0$, such that $f_Y$ satisfies a limit
theorem for the normalization $B_n$: there exists a random
variable $Z$ such that, for every $t\in \R$,
\begin{equation}
  \label{limite_mixing}
  E_Y\left(e^{it \frac{S^Y_{\lfloor n m(Y) \rfloor}  f_Y}{B_n}}\right)
  \to E\left(e^{itZ}\right).
  \end{equation}
\item
\label{hypothese_3}
There exists $b>0$ such that, in the natural extension of $T_Y$,
$\frac{1}{N^b} \sum_0^{N-1} f_Y(T_Y^k y)$ tends almost everywhere
to $0$ when $N \to \pm \infty$.
\item
\label{hypothese_4}
There exists $B'_n=O(B_n^{1/b})$ such that $  \frac{S_n^Y \phi -
nE_Y(\phi)}{B'_n}$ converges in
  distribution.
\end{enumerate}

Then the function $f$ satisfies also a limit theorem:
\begin{equation*}
  E\left(e^{it \frac{S_n f}{B_n}}\right)
  \to E(e^{it Z}),
  \end{equation*}
i.e., $\frac{S_n f}{B_n}$ tends in distribution to $Z$.
\end{thm}

The first assumption is apparently different from the first
assumption in \cite[Theorem A.1]{gouezel:skewproduct}. However,
they are equivalent by \cite{eagleson} (see also
\cite{melbourne_torok}).

\begin{rmq}
\label{rmq_limit_flows} An analogous theorem holds in the
case of flows, when $Y$ is a Poincar\'e section of the flow and
$\phi$ is the return time to this Poincar\'e section, with the same
proof. Since a Poincar\'e section has usually zero measure, it has
to be formulated slightly differently: $E_Y$ will be the
expectation with respect to the probability measure induced by $m$
on $Y$, and in \eqref{limite_mixing} $m(Y)$ should be replaced
with $1/E_Y(\phi)$. Finally, the sums (in the definition of $f_Y$,
and in the definition of the Birkhoff sums of $f$) should be
replaced with integrals, and correspondingly, the normalizing
sequences $B_n$($B_n'$) with appropriate functions $B(T)$
($B'(T)$), $B:\R^+\to \R^+$.
\end{rmq}

\subsection{Proof of Theorem~\ref{main_theorem} for functions
satisfying $(P1)$}

Let $f_0:X_0 \to \R$ be H\"older continuous and satisfy $(P1)$. In
particular, $I\not=0$. Define as in
Section~\ref{section_background} functions $f, \bar f, \bar g$ and
$g$. Since $f$ satisfies \eqref{domaine_attraction} and $\bar g
-\bar f$ is bounded, we obtain $\mu_\Delta( |g|>x) \sim x^{-2}
l(x)$, where
  \begin{equation*}
  l(x)=\frac{I^2 \ell^2}{2\pi}.
  \end{equation*}
By Paragraph~\ref{par_attraction}, the function $g$ is in the
nonstandard domain of attraction of the normal law. More
precisely, set
  \begin{equation*}
  L(x)= \frac{I^2 \ell^2}{\pi} \log(x) \sim
  2\int_1^x \frac{l(u)}{u}.
  \end{equation*}
The functions $l$ and $L$ are the tail functions of $g$,
as defined in Paragraph~\ref{par_attraction}.

Proposition~\ref{estime_integrale_precise} gives
  \begin{equation*}
  \int g (e^{itG}-1)=(y-1) it L(1/|t|) + o(tL(1/|t|)),
  \end{equation*}
where $y=\frac{1}{1-\frac{3}{4}\log 3}$. Moreover, the function
$g$ is locally H\"older on $\Delta$, by \eqref{g_Lipschitz}. Hence, all the
assumptions of Theorem~\ref{thm_ad_young} are satisfied, for
$a=y-1>0$. Let
  \begin{equation*}
  B_n =\sqrt{ n \log n\frac{(2y-1)I^2 \ell^2}{2\pi}},
  \end{equation*}
it satisfies $\frac{n}{B_n^2} (2a+1)L(B_n) \to 1$. Hence, by
Theorem \ref{thm_ad_young}, we obtain that $\frac{\sum_{k=0}^{n-1}
g\circ U^k}{B_n} \to \boN(0,1)$ in distribution with respect to
$\mu_\Delta$. This is equivalent to the same convergence for $\bar
g$, with respect to $\mu_{\bar \Delta}$, since $\bar
g=g\circ\pi_\Delta$ and $\mu_\Delta=(\pi_\Delta)_*(\mu_{\bar
\Delta})$. Since $\bar f$ is cohomologous to $\bar g$, we get the
same convergence for $\bar f$. Finally, since $\bar f=f\circ
\pi_X$ and $\mu=(\pi_X)_*(\mu_{\bar \Delta})$, we get that
  \begin{equation*}
  \frac{\sum_{k=0}^{n-1} f\circ T^k}{B_n} \to \boN(0,1)
  \end{equation*}
on $X$, with respect to $\mu$.

The same argument applies to $\phi_+ - \int \phi_+$, and we get
that $ \frac{\sum_{k=0}^{n-1} \phi_+ \circ T^k - n
\int\phi_+}{B_n}$ converges in distribution.
Hence, Theorem~\ref{thm_probabiliste_general} applies, with $b=1$.

Set $B'_n=B_{ \lfloor n \mu_0(X) \rfloor}$. Since
$\mu_0(X)=\frac{\pi}{2(\pi+\ell)}$ by \eqref{calcule_mu0_X}, we
get
  \begin{equation*}
  B'_n \sim \sqrt{ n \log n \frac{(2y-1)I^2 \ell^2}{4 (\pi+\ell)}}.
  \end{equation*}
Theorem~\ref{thm_probabiliste_general} yields
  \begin{equation*}
  \frac{\sum_{k=0}^{n-1} f_0 \circ T_0^k}{B'_n} \to \boN(0,1).
  \end{equation*}
This concludes the proof of Theorem~\ref{main_theorem}.

\subsection{Proof of Theorem~\ref{main_thm_2}}
\label{sec:gordin}

Let $f_0:X_0 \to \R$ be H\"older continuous with $\int f_0=0$ and $I=0$.
In this case, we can not use the cohomology trick any more, since
the proofs of Lemmas~\ref{u_bornee} and~\ref{lem:g_Holder} relied
heavily on the property $(P1)$. The argument will be to induce on the
basis of the tower $\bar\Delta$, prove a central limit theorem
here (using Gordin's martingale argument), and then get back to
the original space by using Theorem~\ref{thm_probabiliste_general}
twice. The main difference in the inducing process with the
previous paragraph is that we can no more apply
Theorem~\ref{thm_probabiliste_general} with $b=1$. Hence, we will
need to prove that $\frac{1}{|n|^b}\sum_{k=0}^{n-1} f\circ T^k$
converges almost everywhere to $0$, for some $b<1$. Many arguments
of this paragraph are strongly inspired by \cite{lsyoung:annals},
with additional technical complications due to the fact that our
functions are not bounded.

Let $\bar \Delta_0$ be the basis of the tower $\bar \Delta$, and
let $\bar U_0$ be the induced map on $\bar \Delta_0$ (with a
return time $\phi$). Define a new function $\bar f_0$ on $\bar
\Delta_0$, by $\bar f_0(x)=\sum_{k=0}^{\phi(x)-1} \bar f(\bar U^k
x)$.

\begin{lem}
\label{CLT_basis}
There exists $\sigma_0^2\geq 0$ such that
  \begin{equation*}
  \frac{\sum_{k=0}^{n-1} \bar f_0 \circ \bar U_0^k}{\sqrt{n}} \to
  \boN(0,\sigma_0^2).
  \end{equation*}
\end{lem}
\begin{proof}
Since $I=0$, it is not hard to check that there exists
$\alpha_1<1$ such that $|f| \leq n^{\alpha_1}$ on the set of
points bouncing $n$ times along the segments of the stadium. This
implies that there exists $\epsilon_1>0$ such that $f \in
L^{2+\epsilon_1}(X)$. Hence, $\bar f \in L^{2+\epsilon_1}(\bar
\Delta)$. Since the return time $\phi$ belongs to $L^p$ for all
$p<\infty$, we get $\bar f_0 \in L^{2+\epsilon_2}(\bar \Delta_0)$
for some $\epsilon_2>0$.

Let $\Delta_0$ be obtained by identifying the points on the same
stable leaf. It is the basis of the expanding Young tower
$\Delta$. Let $\pi_0:\bar \Delta_0 \to \Delta_0$ be the canonical
projection, and $U_0$ the dynamics induced by $\bar U_0$ on
$\Delta_0$. Let $\boB_0$ be the $\sigma$-algebra on $\bar\Delta_0$
obtained by pulling by $\pi_0$ the $\sigma$-algebra on $\Delta_0$.
A measurable subset $B$ of $\bar\Delta_0$ is $\boB_0$-measurable
if, for almost all $x\in B$, the stable leaf through $x$ is
contained in $B$.

We will prove
  \begin{equation}
  \label{eq:gordin2}
  \sum_{n\geq 0} \norm{E(\bar f_0 \tq \bar U_0^{n}\boB_0)-\bar f_0}_{L^2} <
  \infty
  \end{equation}
and
  \begin{equation}
  \label{eq:gordin1}
  \sum_{n\geq 0} \norm{E(\bar f_0 \tq \bar U_0^{-n}\boB_0)}_{L^2} < \infty.
  \end{equation}
By Gordin's Theorem \cite{gordin}, this will imply the conclusion
of the lemma.

The basis $\bar \Delta_0$ corresponds to a rectangle $R$ for the
dynamics $T$, which is naturally partitioned as $R=\bigcup R_i$,
where $R_i$ is an $s$-subrectangle of $R$. Let $\bar \Delta_{0,i}$
be the corresponding subset of $\bar\Delta_0$, so that $\{
\bar\Delta_{0,i}\}$ gives a partition of $\bar \Delta_0$. Define a
function $A:\bar \Delta_0 \to \R$ by $A(x)=\sum_{k=0}^{\phi(x)-1}
\phi_+( \pi_X \bar U^k x)$. It is constant on each set
$\bar\Delta_{0,i}$, and corresponds to the 
number of times the original map $T_0$ is to be applied to $R_i$ so that 
this s-subrectangle makes a full (Markov) return to the base $R$. 
Since $\phi$ belongs to every $L^p(\bar \Delta_0)$
for $p\geq 1$ and $\phi_+ \in L^p(X)$ for $1\leq p<2$, the function
$A$ belongs to $L^p(\bar \Delta_0)$ for $1\leq p<2$. If $x,y$ are on the same
unstable leaf in a rectangle $\bar \Delta_{0,i}$, we have
  \begin{equation}
  \label{small_unstable}
  |\bar f_0(x) -\bar f_0(y)| \leq C A(x) \tau^{s(x,y)}
  \end{equation}
for some constant $C>0$ and some constant $\tau<1$. Here, $s(x,y)$
is the separation time of $x$ and $y$. Moreover, if $x,y$ are on
the same stable leaf in a rectangle $\bar \Delta_{0,i}$,
  \begin{equation}
  \label{small_stable}
  | \bar f_0(x) -\bar f_0(y)|\leq C A(x) d(\pi_X x, \pi_X y)^\alpha.
  \end{equation}
for some $\alpha>0$.

Since the stable leaves are contracted at each iteration by at
least $\lambda<1$, the atoms of the $\sigma$-algebra $\bar U_0^n
\boB_0$ have a diameter at most $C\lambda^{n}$. By
\eqref{small_stable}, we get
  \begin{equation}
  \label{borne_g2}
  \left| \bar f_0(x)- E(\bar f_0 \tq \bar U_0^n \boB_0)(x) \right|
  \leq C A(x) \lambda^{\alpha n}.
  \end{equation}
Unfortunately, $A$ does not belong to $L^2$, so a further argument
is required to get \eqref{eq:gordin2}. Let $p>0$ be such that
$\frac{1}{p}+\frac{1}{2+\epsilon_2} =\frac{1}{2}$. By
\eqref{borne_g2},
  \begin{equation*}
  1_{A \leq n^{2p}} \left| \bar f_0(x)- E(\bar f_0 \tq \bar U_0^n
  \boB_0)(x) \right| \leq C n^{2p} \lambda^{\alpha n}.
  \end{equation*}
Hence, this series is summable in $L^2$. Moreover,
  \begin{equation*}
  \norm{ 1_{A>n^{2p}} \bar f_0}_{L^2} \leq \norm{1_{A>n^{2p}}}_{L^p}
  \norm{ \bar f_0}_{L^{2+\epsilon_2}}
  \leq \frac{ \left(\int A\right)^{1/p}}{n^2} \norm{ \bar
  f_0}_{L^{2+\epsilon_2}}.
  \end{equation*}
The function $E(\bar f_0 \tq \bar U_0^n \boB_0)$ is bounded in
$L^{2+\epsilon_2}$ by $\norm{\bar f_0}_{L^{2+\epsilon_2}}$. Hence,
we obtain
  \begin{equation*}
  \norm{ 1_{A>n^{2p}} \left| \bar f_0 -E(\bar f_0 \tq \bar U_0^n
  \boB_0) \right| }_{L^2} = O(1/n^2),
  \end{equation*}
which is summable. This proves \eqref{eq:gordin2}.

Let $\bar h=E(\bar f_0 \tq \boB_0)$. This function is constant
along the stable leaves, and has zero integral (since $\bar f_0$
also has zero integral). Hence, it induces a function $h$ on the
quotient $\Delta_0$. Since $\bar f_0\in L^2$, it satisfies $h\in
L^2(\Delta_0)$. The following lemma is an easy consequence of the
H\"older properties of the invariant measure and
\eqref{small_unstable}, see \cite[Sublemma page
612]{lsyoung:annals} for details.
\begin{lem}
There exists constants $C>0$ and $\tau<1$ such that, for all $x,y$
in the same unstable leaf of a set $\bar \Delta_{0,i}$,
  \begin{equation*}
  | \bar h(x) - \bar h(y)| \leq C A(x) \tau^{s(x.y)}.
  \end{equation*}
\end{lem}
The function $A$ is integrable. Hence, by \cite[Lemma
3.4]{gouezel:stable}, this implies that the function $\hat{U}_0 h$
is H\"older continuous on $\Delta_0$. By \cite[Corollary
3.3]{gouezel:stable}, we get:
  \begin{equation}
  \label{dec_exp}
  \hat{U}_0^n h \text{ tends exponentially fast to $0$
  in the space of H\"older continuous functions on $\Delta_0$.}
  \end{equation}
A computation gives
  \begin{equation*}
  \norm{ E(\bar f_0 \tq \bar U_0^{-n} \boB_0)}_{L^2}^2
  = \int h \cdot (\hat{U}_0^n h)\circ U_0^n
  \leq \norm{h}_{L^2} \bigl\| (\hat{U}_0^n h)\circ U_0^n \bigr\|_{L^2}
  =  \norm{h}_{L^2} \bigl\| \hat{U}_0^n h\bigr\|_{L^2}.
  \end{equation*}
Hence, this term is exponentially small. This proves
\eqref{eq:gordin1} and concludes the proof of
Lemma~\ref{CLT_basis}.
\end{proof}

The return time $\phi$ also satisfies a central limit theorem, by
the same argument. Hence, by
Theorem~\ref{thm_probabiliste_general} (applied with $b=1$), there
exists $\sigma_1^2 \geq 0$ such that
  \begin{equation*}
  \frac{ \sum_{k=0}^{n-1} \bar f \circ \bar U^k}{\sqrt{n}} \to
  \boN(0,\sigma_1^2).
  \end{equation*}
Going from $\bar \Delta$ to $X$, it implies that
  \begin{equation*}
  \frac {\sum_{k=0}^{n-1} f\circ T^k}{\sqrt{n}} \to
  \boN(0,\sigma_1^2).
  \end{equation*}
Moreover, the return time $\phi_+ : X \to \N$ satisfies a limit
theorem with normalization $\sqrt{n\log n}$. Since $\sqrt{n}=o
(\sqrt{n \log n})$, we can unfortunately not apply Theorem
\ref{thm_probabiliste_general} with $b=1$. However, if we can
prove the following lemma, then this theorem applies with $b<1$,
and this concludes the proof of Theorem~\ref{main_thm_2}.

\begin{lem}
\label{lem_strong_ergthm_f}
For all $b>1/2$,
  \begin{equation}
  \label{conv}
  \frac{1}{|n|^b}\sum_{k=0}^{n-1} f\circ T^k \to 0
  \end{equation}
almost everywhere in $X$ when $n \to \pm \infty$.
\end{lem}
\begin{proof}
We first estimate the decay of correlations of $\bar f_0$  for
$\bar U_0$. We will use the notations of the proof of Lemma
\ref{CLT_basis}. We have
  \begin{equation}
  \label{decay_cor}
  \int \bar f_0 \cdot \bar f_0 \circ \bar U_0^{2n}
  =\int \bar f_0 \cdot E(\bar f_0\circ \bar U_0^n \tq \boB_0)\circ \bar
  U_0^n
  + \int \bar f_0 \cdot 
  \left(  \bar f_0 \circ \bar U_0^{2n} - E( \bar f_0\circ
  \bar U_0^n\tq \boB_0)\circ \bar U_0^n \right).
  \end{equation}
The contraction properties of $\bar U_0$ along stable manifolds
and \eqref{small_stable} give $|\bar f_0 \circ \bar U_0^n (x)-E( \bar
f_0\circ \bar U_0^n \tq \boB_0)(x)|\leq C A(\bar U_0^n x)
\lambda^{\alpha n}$. Hence, the second integral in
\eqref{decay_cor} is at most
  \begin{equation*}
  \int |\bar f_0| \cdot A\circ \bar U_0^{2n} \lambda^{\alpha n}
  \leq \norm{\bar f_0}_{L^{2+\epsilon_2}} \norm{A}_{L^p}
  \lambda^{\alpha n},
  \end{equation*}
where $p<2$ is chosen so that
$\frac{1}{2+\epsilon_2}+\frac{1}{p}=1$. Hence, this term decays
exponentially fast.

In the first integral of \eqref{decay_cor}, the function $E( \bar
f_0\circ \bar U_0^n \tq \boB_0)\circ \bar U_0^n$ is
$\boB_0$-measurable (i.e., constant along stable leaves). Hence,
this integral is equal to
  \begin{equation}
  \label{dec_cor2}
  \int \bar h \cdot E(\bar f_0\circ \bar U_0^n \tq \boB_0)\circ \bar
  U_0^n.
  \end{equation}
Let $\bar h_n= E( \bar f_0\circ \bar U_0^n\tq \boB_0)$, it is
$\boB_0$-measurable and defines a function $h_n$ on the quotient
$\Delta_0$. The integral \eqref{dec_cor2} is then equal to
  \begin{equation}
  \label{dec_cor3}
  \int_{\Delta_0} h \cdot h_n \circ U_0^n
  = \int \hat{U}_0^n h \cdot h_n.
  \end{equation}
The $L^2$-norm of $h_n$ is bounded independently of $n$. By
\eqref{dec_exp}, \eqref{dec_cor3} is exponentially small. This
proves that $\int \bar f_0 \cdot \bar f_0\circ \bar U_0^{2n}$
decays exponentially. In the same way, $\int \bar f_0 \cdot \bar
f_0\circ \bar U_0^{2n+1}$ decays exponentially.

Since the correlations of $\bar f_0$ decay exponentially fast and
$\bar f_0\in L^2$,
 \cite[Theorem 16]{vitesse_birkhoff} implies that
$\frac{1}{n^b}\sum_{k=0}^{n-1} \bar f_0 \circ \bar U_0^k$ tends to
zero almost everywhere when $n\to +\infty$, for all $b>1/2$.

Now to see that $\frac{1}{n^b} \sum_{k=0}^{n-1} \bar f\circ \bar
U^k $ tends to zero almost everywhere in $\bar \Delta$  when $n\to
+\infty$, for all $b>1/2$, we use \cite[Lemma 2.1
(a)]{melbourne_torok} which gives this convergence on $\bar
\Delta_0$. However, by the ergodicity of $\bar U$, the set on which
this convergence holds must have either full or zero measure. As
$\bar \Delta_0$ has positive measure, we get this convergence
almost everywhere on $\bar \Delta$. Finally, this implies the same for
$f$ in $X$. We have proved \eqref{conv} for any $b>1/2$ when $n\to
+\infty$.

To deal with $n \to -\infty$, we go to the natural extension. It
is sufficient to prove the result for $\bar f_0$ in $\bar
\Delta_0$, since the previous reasoning still applies (using the
fact that the natural extension is functorial, i.e., the natural
extension commutes with induction and projections). In the natural
extension $\bar \Delta'_0$ of $\bar \Delta_0$, we have $\int \bar
f'_0 \cdot \bar f'_0 \circ {{\bar U}'_0}{}^{-n} = \int \bar f_0
\circ \bar U_0^n \cdot \bar f_0$, which is exponentially small.
Hence, \cite[Theorem 16]{vitesse_birkhoff} still applies and gives
the desired result.
\begin{comment}
We could also say that $T^{-1}$ has the same properties as $T$.
Hence, the same arguments as above apply to $T^{-1}$.
\end{comment}
\end{proof}

\begin{rmq}
\label{rmq:strong_ergthm_f0} As $\mu_0(X)>0$, we may apply
\cite[Lemma 2.1 (a)]{melbourne_torok} just as we did in the proof
above to see that Lemma~\ref{lem_strong_ergthm_f} implies
  \begin{equation*}
  \frac{1}{|n|^b}\sum_{k=0}^{n-1} f_0\circ T_0^k \to 0
  \end{equation*}
almost everywhere when $n \to \pm \infty$, for any $b>1/2$.
\end{rmq}

\section*{Acknowledgements}

We are much grateful to D. Sz\'asz and T. Varj\'u for useful
discussions and for their valuable remarks on earlier versions of
the manuscript. This paper has grown out of discussions we had at
the CIRM conference on multi-dimensional non-uniformly hyperbolic
systems in Marseille in May 2004, and while S.\ G.\ visited the
Institute of Mathematics of the BUTE in October 2004. The
hospitality of both institutions, along with the financial support
of Hungarian National Foundation for Scientific Research (OTKA),
grants T32022 and TS040719 is thankfully acknowledged.

\appendix
\section{Proof of Lemma~\ref{renouvell_exp}}

Let $U_0$ be the map induced by $U$ on the basis $\Delta_0$ of the
tower. Denote by $\phi$ the first return time on the basis, so
that $U_0(x)=U^{\phi(x)}(x)$. Note that $\phi(x)$ can also be
defined for $x\in \Delta\setminus \Delta_0$ as the first hitting
time of the basis.

Let $F$ be a finite subset of $\N$. Let $(n_i)_{i\in F}$ be
positive integers. Let
  \begin{equation*}
  K(F,n_i)=\{x \in \Delta_0 \tq \forall i \in F, \phi(U_0^i x)=n_i\}.
  \end{equation*}
\begin{lem}
\label{lem_KFn}
There exists a constant $C$ such that, for all $F$ and $n_i$ as
above,
  \begin{equation*}
  \mu_\Delta( K(F,n_i)) \leq \prod_{i\in F} (C \rho^{n_i}).
  \end{equation*}
\end{lem}
\begin{proof}
The proof is by induction on $\max F$, and the result is trivial
when $F=\emptyset$.  Write $F'=\{ i-1 \tq i\in F, i\geq 1\}$ and,
for $i\in F'$, set $n'_i=n_{i+1}$.

if $0\not \in F$, $K(F, n_i)= U_0^{-1}(K(F', n_i'))$. Since $U_0$
preserves $\mu_{\Delta}$ and $\max F' <\max F$, we get the result.
Otherwise, $0\in F$. Then $K(F,n_i)=U_0^{-1}(K(F',n'_i))\cap \{
x\in \Delta_0, \phi(x)=n_0\}$. By bounded distortion, we get
  \begin{equation*}
  \mu_\Delta (K(F,n_i)) \leq C \mu_\Delta( K(F',n_i'))
  \mu_\Delta \{x\in \Delta_0, \phi(x)=n_0\}
  \leq C \mu_\Delta(K(F',n'_i)) \rho^{n_0}.
  \qedhere
  \end{equation*}
\end{proof}

\begin{lem}
\label{lem_renouv_temp}
There exist $C>0$ and $\theta<1$ such that, for all $n \in \N$,
  \begin{equation*}
  \int_{U^{-n} \Delta_0} \tau^{\Psi_n} \leq C \theta^n.
  \end{equation*}
\end{lem}
\begin{proof}
Let $\kappa>0$ be very small (how small will be specified later in
the proof). Then
  \begin{equation*}
  U^{-n} \Delta_0 \subset \{ x \in \Delta \tq \Psi_n(x)\geq \kappa n\}
  \cup \{x \in \Delta \tq \phi(x) \geq n/2\}
  \cup \{ x\in \Delta \tq \phi(x)<n/2, \Psi_n(x) < \kappa n\}.
  \end{equation*}
On the first of these sets, $\tau^{\Psi_n} \leq \tau^{\kappa n}$,
whence the integral of $\tau^{\Psi_n}$ is exponentially small. The
second of these sets has exponentially small measure. Finally, the
last of these sets is contained in $\bigcup_{i=0}^{n/2} U^{-i}
\Gamma_n$, where
  \begin{equation*}
  \Gamma_n =\{ x\in \Delta_0 \tq \sum_{0\leq i \leq \kappa n}\phi(U_0^i
  x) \geq n/2\}.
  \end{equation*}
To conclude the proof of the Lemma, it is sufficient to prove that
the measure of $\Gamma_n$ is exponentially small.

Take $L\in \N$ such that $\forall n\geq L, (C \rho)^n \leq
\rho^{n/2}$, where $C$ is the constant given by
Lemma~\ref{lem_KFn}. For $x\in \Gamma_n$, let $F(x):=\{ 0\leq i
\leq \kappa n \tq \phi(U_0^i x) \geq L\}$. Then
  \begin{equation*}
  \sum_{i\in F(x)} \phi(U_0^i x) \geq \frac{n}{2} -\sum_{i\not \in
  F(x)} L \geq (1/2 -L \kappa) n.
  \end{equation*}
This implies that
  \begin{equation*}
  \Gamma_n \subset \bigcup_{ F \subset [0, \lfloor \kappa n \rfloor]}
  \bigcup_{ \substack{n_i \geq L \\ \sum_{i\in F} n_i \geq
  (1/2-L\kappa)n}} K(F,n_i).
  \end{equation*}
By Lemma~\ref{lem_KFn}, we get
  \begin{align*}
  \mu_\Delta(\Gamma_n) &\leq \sum_{F \subset [0, \lfloor \kappa n \rfloor]}
  \sum_{ \substack{n_i \geq L \\ \sum_{i\in F} n_i \geq
  (1/2-L\kappa)n}} \prod_{i \in F} (C \rho^{n_i})
  \leq \sum_{k=0}^{\lfloor \kappa n \rfloor} \binom{ \lfloor \kappa n
  \rfloor}{k} \sum_{ \substack{ n_0,\dots,n_{k-1} \geq L \\ \sum n_i
  \geq (1/2-L\kappa)n}} (C \rho^{n_0})\dots (C\rho^{n_{k-1}})
  \\&
  \leq 2^{\kappa n} \sum_{0\leq k \leq \kappa n}
  \sum_{ \substack{ n_0,\dots,n_{k-1} \geq L \\ \sum n_i
  \geq (1/2-L\kappa)n}} \rho^{\sum n_i /2}
  \leq 2^{\kappa n} \sum_{0\leq k \leq \kappa n}
  \sum_{ \substack{n_0,\dots,n_{k-1} \in \N \\
  \sum n_i \geq (1/2-L\kappa)n}} \rho^{\sum n_i /2}.
  \end{align*}
For $r\in\N$,
  \begin{equation*}
  \sum_{n_0+\dots+n_{k-1}=r} \rho^{ \sum n_i /2}= \rho^{r/2}
  \Card\{ n_0,\dots,n_{k-1} \tq \sum n_i =r\} = \rho^{r/2}
  \binom{r+k}{k} \leq \rho^{r/2} \frac{(r+k)^k}{k!}.
  \end{equation*}
Hence,
  \begin{equation*}
  \mu_\Delta(\Gamma_n)
  \leq 2^{\kappa n} \sum_{0\leq k \leq \kappa n} \sum_{r\geq
  (1/2-L\kappa)n} \rho^{r/2} \frac{(r+k)^k}{k!}.
  \end{equation*}
The sequence $u_r= \rho^{r/2} \frac{(r+k)^k}{k!}$ satisfies
$\frac{u_{r+1}}{u_r} \leq \rho':= \rho^{1/2}
e^{\frac{\kappa}{1/2-L\kappa}}$ for all $r \geq (1/2-L\kappa)n$
and $k \leq \kappa n$. if $\kappa$ is small enough, $\rho'<1$, and
we get
  \begin{equation*}
  \mu_\Delta(\Gamma_n) \leq 2^{\kappa n} \sum_{0\leq k \leq \kappa n} \rho^{
  (1/2-L\kappa) n /2} \frac{ \bigl((1/2-L\kappa)n + \kappa n\bigr)^k}
  {k!}\frac{1}{1-\rho'}
  \leq \frac{2^{\kappa n}}{1-\rho'} \rho^{
  (1/2-L\kappa) n /2}  \sum_{0\leq k \leq \kappa n} \frac{n^k}{k!}.
  \end{equation*}
The sequence $\frac{n^k}{k!}$ is increasing for $k\leq n$. Hence,
we finally get
  \begin{equation*}
  \mu_\Delta(\Gamma_n) \leq \frac{2^{\kappa n}}{1-\rho'} \rho^{
  (1/2-L\kappa) n /2} (\kappa n+1) \frac{n^{\lfloor \kappa n\rfloor}}{
  \lfloor \kappa n \rfloor !}.
  \end{equation*}
Using Stirling's Formula, it is easy to check that this expression
is exponentially small if $\kappa$ is small enough. This concludes
the proof.
\end{proof}

\begin{proof}[Proof of Lemma~\ref{renouvell_exp}]
Let $\theta$ be given by Lemma~\ref{lem_renouv_temp}. Choose
$\alpha>0$ so that $e^{\epsilon \alpha} \theta <1$. Then
  \begin{equation*}
  U^{-n} \Delta_0 \subset \{ x\in \Delta \tq \haut(x) \geq \alpha n\}
  \cup \Bigl[ \{ x \in \Delta \tq \haut(x) < \alpha n\}  \cap
  U^{-n}\Delta_0 \Bigr].
  \end{equation*}
Hence,
  \begin{equation*}
  \int_{U^{-n} \Delta_0} e^{\epsilon \omega} \tau^{\Psi_n}
  \leq \int_{ \omega\geq \alpha n} e^{\epsilon \omega}
  + e^{\epsilon \alpha n} \int_{U^{-n} \Delta_0} \tau^{\Psi_n}.
  \end{equation*}
The first term is exponentially small since $e^{\epsilon}\rho<1$.
Lemma~\ref{lem_renouv_temp} and the definition of $\alpha$ also
imply that the second term is exponentially small.
\end{proof}

\bibliography{biblio}

\begin{thebibliography}{SV04b}

\bibitem[Aar97]{aaronson:book}
Jon Aaronson.
\newblock {\em An introduction to infinite ergodic theory}, volume~50 of {\em
  Mathematical Surveys and Monographs}.
\newblock American Mathematical Society, 1997.

\bibitem[AD01]{aaronson_denker:central}
Jon Aaronson and Manfred Denker.
\newblock {A local limit theorem for stationary processes in the domain of
  attraction of a normal distribution.}
\newblock In N.~Balakrishnan, I.A. Ibragimov, and V.B. Nevzorov, editors, {\em
  Asymptotic methods in probability and statistics with applications. Papers
  from the international conference, St. Petersburg, Russia, 1998}, pages
  215--224. Birkh{\"a}user, 2001.

\bibitem[Bun79]{bunimovich:stadium}
Leonid Bunimovich.
\newblock On the ergodic properties of nowhere dispersing billiards.
\newblock {\em Comm. Math. Phys.}, 65:295--312, 1979.

\bibitem[BY93]{baladi_young}
Viviane Baladi and Lai-Sang Young.
\newblock On the spectra of randomly perturbed expanding maps.
\newblock {\em Comm. Math. Phys.}, 156:355--385, 1993.

\bibitem[Che97]{chernov:entropy}
Nikolai Chernov.
\newblock Entropy, {L}yapunov exponents and mean free path for billiards.
\newblock {\em J. Statist. Phys.}, 88:1--29, 1997.

\bibitem[Che99]{chernov:decay}
Nikolai Chernov.
\newblock Decay of correlations and dispersing billiards.
\newblock {\em J. Statist. Phys.}, 94:513--556, 1999.

\bibitem[CZ]{chernov:slow}
Nikolai Chernov and Hongkun Zhang.
\newblock Billiards with polynomial mixing rates.
\newblock Preprint.

\bibitem[Eag76]{eagleson}
G.~K. Eagleson.
\newblock Some simple conditions for limit theorems to be mixing.
\newblock {\em Teor. Verojatnost. i Primenen.}, 21(3):653--660, 1976.

\bibitem[Gor69]{gordin}
Mikhail Gordin.
\newblock The central limit theorem for stationary processes.
\newblock {\em Dokl. Akad. Nauk SSSR}, 188:739--741, 1969.

\bibitem[Gou03]{gouezel:skewproduct}
S{\'e}bastien Gou{\"e}zel.
\newblock Statistical properties of a skew-product with a curve of neutral
  points.
\newblock Preprint, 2003.

\bibitem[Gou04]{gouezel:stable}
S{\'e}bastien Gou{\"e}zel.
\newblock Central limit theorem and stable laws for intermittent maps.
\newblock {\em Probab. Theory and Rel. Fields}, 128:82--122, 2004.

\bibitem[Hen93]{hennion}
Hubert Hennion.
\newblock Sur un th{\'e}or{\`e}me spectral et son application aux noyaux
  lipschitziens.
\newblock {\em Proc. Amer. Math. Soc.}, 118:627--634, 1993.

\bibitem[Kac96]{vitesse_birkhoff}
A.~G. Kachurovski\u\i.
\newblock Rates of convergence in ergodic theorems.
\newblock {\em Russian Math. Surveys}, 51:653--703, 1996.

\bibitem[KL99]{keller_liverani}
Gerhard Keller and Carlangelo Liverani.
\newblock Stability of the spectrum for transfer operators.
\newblock {\em Ann. Scuola Norm. Sup. Pisa Cl. Sci. (4)}, 28(1):141--152, 1999.

\bibitem[Mac83]{machta:cusp}
J.~Machta.
\newblock Power law decay of correlations in a billiard problem.
\newblock {\em J. Statist. Phys.}, 32:555--564, 1983.

\bibitem[Mar04]{markarian:slow}
Roberto Markarian.
\newblock Billiards with polynomial decay of correlations.
\newblock {\em Ergodic Theory Dynam. Systems}, 24:177--197, 2004.

\bibitem[MT04]{melbourne_torok}
Ian Melbourne and Andrew T{\"o}r{\"o}k.
\newblock Statistical limit theorems for suspension flows.
\newblock {\em Israel J. Math.}, pages 191--210, 2004.

\bibitem[SV]{szasz_varju:infinite}
Domokos {Sz}\'asz and Tam\'as Varj\'u.
\newblock in preparation.

\bibitem[SV04a]{szasz_varju:finite}
Domokos {Sz}\'asz and Tam\'as Varj\'u.
\newblock Local limit theorem for the {L}orentz process and its recurrence on
  the plane.
\newblock {\em Ergodic Theory Dynam. Systems}, 24:257--278, 2004.

\bibitem[SV04b]{szasz_varju:expository}
Domokos {Sz}\'asz and Tam\'as Varj\'u.
\newblock Markov towers and stochastic properties of billiards.
\newblock In {\em Modern dynamical systems and applications}, pages 433--445.
  Cambridge University Press, 2004.

\bibitem[You98]{lsyoung:annals}
Lai-Sang Young.
\newblock Statistical properties of dynamical systems with some hyperbolicity.
\newblock {\em Ann. of Math. (2)}, 147:585--650, 1998.

\bibitem[You99]{lsyoung:recurrence}
Lai-Sang Young.
\newblock Recurrence times and rates of mixing.
\newblock {\em Israel J. Math.}, 110:153--188, 1999.

\end{thebibliography}
\bibliographystyle{alpha}

\end{document}